\documentclass[oneside,12pt]{amsart}
\usepackage{amssymb, mathrsfs, amscd}
\usepackage[all]{xy}

\usepackage{graphics}
\usepackage{graphicx}
\usepackage{enumerate}

\newtheorem*{thma}{Theorem A}
\newtheorem*{thmb}{Theorem B}
\newtheorem*{thmc}{Theorem C}

\newtheorem{theorem}{Theorem}[section]
\newtheorem{definition}[theorem]{Definition}
\newtheorem{lemma}[theorem]{Lemma}

\newtheorem{proposition}[theorem]{Proposition}
\newtheorem{corollary}[theorem]{Corollary}
\newtheorem*{corollary*}{Corollary}

\newtheorem{remark}[theorem]{Remark}

\newtheorem{observation}[theorem]{Observation}

\newtheorem{question}[theorem]{Question}

\author{Dieter Degrijse}
\address{Department of Mathematics, KU Leuven, Kortrijk, Belgium}%
\email{Dieter.Degrijse@kuleuven-kulak.be}%

\title[]{Bredon cohomological dimensions for proper actions and Mackey functors}
\date{\today}

\newcommand{\mF}{\mathcal {F}}

\newcommand{\orb}{\mathcal{O}_{\mF}G}
\newcommand{\orbmod}{\mbox{Mod-}\mathcal{O}_{\mF}G}

\bigskip
\begin{document}

\maketitle
\begin{abstract}
For groups with a uniform bound on the length of chains of finite subgroups, we study the relationship between the Bredon cohomological dimension for proper actions and the notions of cohomological dimension one obtains by restricting the coefficients of Bredon cohomology to (cohomological) Mackey functors or fixed point functors. We also investigate the closure properties of the class of groups with finite Bredon cohomological dimension for Mackey functors.
\end{abstract}
\section{Introduction}

In this paper we build upon the work of Mart{\'\i}nez-P{\'e}rez and Nucinkis in \cite{MartinezNucinkis06} and investigate Bredon cohomology with Mackey functor coefficients and the notions of cohomological dimension that can be derived from this. Throughout, let $G$ be a discrete group and let $\mathcal{F}$ be the family of finite subgroups of $G$. The Bredon cohomology  $\mathrm{H}^{\ast}_{\mathcal{F}}(G,-)$ of $G$ for the family $\mathcal{F}$ is a generalization of the ordinary group cohomology of $G$, where the category of $G$-modules is replaced by the category $\mbox{Mod-}\mathcal{O}_{\mathcal{F}}G$ of contravariant functors from the orbit category $\mathcal{O}_{\mathcal{F}}G$ to the category of abelian groups. The coefficient category $\mbox{Mod-}\mathcal{O}_{\mathcal{F}}G$ has an interesting subclass (see Remark \ref{remark: subclass}) $ \mathrm{Mack}_{\mathcal{F}}(G) $ consisting of Mackey functors. Bredon cohomology with Mackey functor coefficients is often easier to handle than Bredon cohomology with general coefficients since it has more properties in common with ordinary group cohomology, a key one being the existence of transfer maps. Inside the category $\mathrm{Mack}_{\mathcal{F}}(G)$ one has the subcategory of cohomological Mackey functors $\mathrm{coMack}_{\mathcal{F}}(G)$, of which the category of fixed point functors $\mathcal{F}\mathrm{ix}(G)$ is an important subcategory. One therefore has the following chain of different classes of coefficients for Bredon cohomology
\[   \mathcal{F}\mathrm{ix}(G) \subseteq \mathrm{coMack}_{\mathcal{F}}(G) \subseteq \mathrm{Mack}_{\mathcal{F}}(G) \subseteq \mbox{Mod-}\mathcal{O}_{\mathcal{F}}G.      \]
By considering Bredon cohomology with coefficients in these different classes, one obtains four notions of cohomological dimension, satisfying
\begin{equation}\label{eq: all the inv}\mathcal{F}\mathrm{cd}(G)\leq \underline{\mathrm{cd}}_{co\mathcal{M}}(G) \leq \underline{\mathrm{cd}}_{\mathcal{M}}(G) \leq \underline{\mathrm{cd}}(G).\end{equation}
Here $\underline{\mathrm{cd}}(G)$ is the Bredon cohomological dimension of $G$ for proper actions, and $\mathcal{F}\mathrm{cd}(G)$ is an invariant introduced by Nucinkis in \cite{nucinkis99} (see also \cite{nucinkis00}) using slightly different methods. 

The invariant $\underline{\mathrm{cd}}(G)$ has a geometric interpretation. A proper $G$-CW complex is a $G$-CW-complex with finite isotropy groups. A classifying space for proper actions of $G$, or model for $\underline{E}G$, is a proper $G$-CW-complex $X$ such that $X^H$ is contractible for every finite subgroup $H$ of $G$. These spaces appear naturally throughout geometric group theory, and are important tools for studying groups. Classifying spaces for proper actions also have important K-theoretical applications. They appear for example in the Baum-Connes conjecture  (see e.g. \cite{BCH}, \cite{MV},\cite{LuckReich}) and in a generalization to infinite groups by L\"{u}ck an Oliver of the Atiyah-Segal completion theorem (see \cite{LuckOliver01}). We refer the reader to the survey paper \cite{Luck2} for more information about classifying spaces for families of subgroups and their applications. With these applications and others in mind, it is important to have models for $\underline{E}G$ with good geometric finiteness properties. One such finiteness property is being finite dimensional. The geometric dimension for proper actions of $G$ is by definition the number 
\[\underline{\mathrm{gd}}(G)=\inf\{ \dim(X) \ | \ X \ \mbox{is a model for $\underline{E}G$}\}.\]  This geometric invariant is connected to its algebraic counterpart by the following inequalities (see \cite[0.1]{LuckMeintrup}) \[\underline{\mathrm{cd}}(G) \leq  \underline{\mathrm{gd}}(G) \leq \max\{\underline{\mathrm{cd}}(G) ,3\}.\] We will investigate the relationship between the invariants in (\ref{eq: all the inv}), and study the behavior of these invariants under group extensions and other operations. Our reasons for doing this are twofold.

The first reason is inspired by the fact that one would like an as easy as possible geometric or algebraic criterion that decides whether or not a given group $G$ admits a finite dimensional model for $\underline{E}G$. Let $n_G$ denote the smallest possible dimension of a contractible proper $G$-CW-complex.  A conjecture of Kropholler and Mislin (e.g. see \cite[Conj. 43.1]{GuidoConjBook}) says that $n_G < \infty$ if and only if $\underline{\mathrm{gd}}(G)<\infty$. Kropholler and Mislin verified their conjecture for groups of type $FP_{\infty}$ in \cite[Th. A]{KrophollerMislin}, and later in \cite[Th. 1.10]{Luck1} L\"{u}ck proved it more generally for groups of finite length, i.e.~with a uniform bound on the length of chains of finite subgroups. The general case however, remains unsolved. Nucinkis formulated an algebraic version of this conjecture (see \cite{nucinkis00}) stating that $\mathcal{F}\mathrm{cd}(G)<\infty$ if and only if $\underline{\mathrm{cd}}(G)<\infty$. So, $\mathcal{F}\mathrm{cd}(G)$ should be seen as an algebraic version of $n_G$, although the connection between the two is presently only known to be one-sided. In general one has $\mathcal{F}\mathrm{cd}(G)\leq n_G$, by a result of Kropholler and Wall (see \cite[Th. 1.4]{KrophollerWall11}) but it is not even known whether the finiteness of $\mathcal{F}\mathrm{cd}(G)$ implies the finiteness of $n_G$, for all groups $G$. The invariants $\underline{\mathrm{cd}}_{co\mathcal{M}}(G)$ and $\underline{\mathrm{cd}}_{\mathcal{M}}(G) $ in (\ref{eq: all the inv}) lead to refinements of Nucinkis' conjecture. One could hope to prove the conjecture by using $\underline{\mathrm{cd}}_{co\mathcal{M}}(G) $ and $\underline{\mathrm{cd}}_{\mathcal{M}}(G) $ to interpolate between $\mathcal{F}\mathrm{cd}(G)$ and $\underline{\mathrm{cd}}(G)$, or disprove it by showing that the conjectured implication already fails at $\underline{\mathrm{cd}}_{co\mathcal{M}}(G) $ or $\underline{\mathrm{cd}}_{\mathcal{M}}(G) $. 

Our second reason for studying Bredon (co)homology with Mackey functor coefficients and the accompanying notions of (co)homological dimension is that Bredon (co)homology coefficients appearing in applications often turn out to te Mackey functors. For example, the Baum-Connes conjecture or the Atiyah-Segal completion theorem mentioned above contain generalized $G$-(co)homology theories wich can be computed via an Atiyah-Hirzeburch spectral sequence (see \cite[Th. 4.7]{DL}) whose second page contains Bredon (co)homology groups with Mackey functor coefficients. It is therefore of interest to conduct a general study of Bredon (co)homology with Mackey functor coefficients and the related notions of dimension.\\

The length $l(H)$ of a finite subgroup $H$ of $G$ is the length of the longest chain of subgroups of $H$. The length $l(G)$ of $G$ is the supremum of the lengths of the finite subgroups of $G$. Given $H \in \mathcal{F}$, let $N_G(H)$ be the normalizer of $H$ in $G$. It is known that the Weyl-group $W_G(H)=N_G(H)/H$ satisfies $\mathcal{F}\mathrm{cd}(W_G(H))\leq \mathcal{F}\mathrm{cd}(G)$ (see \cite[Lemma 4.3]{nucinkis00}). Regarding Nuckinkis' conjecture and the relationship between the invariants in $(\ref{eq: all the inv})$, it is known that $\underline{\mathrm{cd}}(G)\leq \mathcal{F}\mathrm{cd}(G)+l(G)$ when $l(G)<\infty$, and $\mathrm{vcd}(G)=\mathcal{F}\mathrm{cd}(G)=\underline{\mathrm{cd}}_{\mathcal{M}}(G)$ when $G$ is virtually torsion-free. This first statement follows from combining Theorem 4.4 of \cite{nucinkis00}, Theorem 1.10 of \cite{Luck1} and Theorem 3.10 of \cite{Martinez07}, and the second one is the main theorem of \cite{MartinezNucinkis06}. In Theorem A, B and C below, we will generalize or improve these results. Recently,  Nucinkis' conjecture for groups with unbounded torsion was studied in \cite{GandiniNucinkis}, and in \cite{Martinez11} and \cite{Martinez12} the invariant $\underline{\mathrm{cd}}(G)$ was investigated using posets of finite subgroups. 
\begin{thma} Let $G$ be a group with $l(G)<\infty$. Then we have

\[     \mathcal{F}\mathrm{cd}(G) = \underline{\mathrm{cd}}_{co\mathcal{M}}(G) = \underline{\mathrm{cd}}_{\mathcal{M}}(G)  \leq \underline{\mathrm{cd}}(G) \leq \max_{H \in \mathcal{F}}\Big\{  \mathcal{F}\mathrm{cd}(W_G(H))+l(H)\Big\}.  \]

\end{thma}

\noindent Regarding these inequalities, it should be mentioned that (see Remark \ref{remark: attained upper bounds}) for each $s \in \mathbb{N}$, there exists a virtually torsion-free group $G_s$ with $l(G_s)=s$ such that 
\[\mathcal{F}\mathrm{cd}(G_s)+s=\underline{\mathrm{cd}}(G_s).\]
Also, in Proposition \ref{prop: Fcd2} we prove in general that \[\mathcal{F}\mathrm{cd}(G) = \underline{\mathrm{cd}}_{co\mathcal{M}}(G)\ \mbox{ if $\underline{\mathrm{cd}}_{co\mathcal{M}}(G)<\infty$}.\]

It is a standard fact that the ordinary cohomological dimension of groups behaves sub-additively, with respect to group extensions. Since there are examples of virtually torsion-free groups $G$ for which $\mathrm{vcd}(G)<\underline{\mathrm{cd}}(G)$ (see \cite[Th. 6 and Ex. 12]{LearyNucinkis}), one concludes that this is no longer true for Bredon cohomological dimension. Moreover, it is not even known whether the class of groups with finite Bredon cohomological dimension is closed under group extensions. Since Bredon cohomology with Mackey functor coefficients behaves better than Bredon cohomology with general coefficients, one could hope that the invariants $\underline{\mathrm{cd}}_{(co)\mathcal{M}}(G)$ (i.e.~ $\underline{\mathrm{cd}}_{co\mathcal{M}}(G)$ and $\underline{\mathrm{cd}}_{\mathcal{M}}(G)$) do behave sub-additively, with respect to group extensions. In the second main result, we prove this under some additional assumptions.
\begin{thmb} Consider a short exact sequence
\[ 1 \rightarrow N \rightarrow G \xrightarrow{\pi} Q \rightarrow 1 \]
such that every finite index overgroup of $N$ in $G$ has a uniform bound on the lengths of its finite subgroups that are not contained in $N$. If either  $ \underline{\mathrm{cd}}_{(co)\mathcal{M}}(G)<\infty$ or $\underline{\mathrm{cd}}(N)<\infty$, then we have 
\[   \underline{\mathrm{cd}}_{(co)\mathcal{M}}(G) \leq   \underline{\mathrm{cd}}_{(co)\mathcal{M}}(N) + \underline{\mathrm{cd}}_{(co)\mathcal{M}}(Q).        \] 
\end{thmb}

Given a class of groups $\mathfrak{C}$, let $\mathrm{H}_1(n)\mathfrak{C}$ be the class of all groups $G$ for which there exists a contractible $G$-CW-complex of dimension at most $n$ with isotropy in $\mathfrak{C}$. As a variation on this, define $\widetilde{\mathrm{H}}_1(n)\mathfrak{C}$ to be the class of all groups $G$ that act by isometries, discrete orbits and with isotropy in $\mathfrak{C}$ on separable CAT($0$)-space of topological dimension at most $n$. Given two classes of groups $\mathfrak{C}$ and $\mathfrak{D}$, let $\mathrm{Ext}(\mathfrak{C},\mathfrak{D})$ the class of all groups that arise as an extension of a group in $\mathfrak{C}$ by a group in $\mathfrak{D}$. Denote by $\mathfrak{F}_d$, respectively $\mathfrak{M}_d$, the class of all groups $G$ with $\mathcal{F}\mathrm{cd}(G) \leq d$, respectively $\underline{\mathrm{cd}}_{\mathcal{M}}(G) \leq d$. Finally, denote by $\mathfrak{L}$ the class of all groups that have finite length. 
By Theorem A, we have $\mathfrak{L}\cap \mathfrak{F}_d= \mathfrak{L}\cap \mathfrak{M}_d$ for all $d \in \mathbb{N}$. We prove that the classes of groups 
\[  \mathfrak{F}= \bigcup_{d \in \mathbb{N}} \mathfrak{F}_d \ \ \mbox{and} \ \ \mathfrak{M} = \bigcup_{d \in \mathbb{N}}  \mathfrak{M}_d\]
satisfy the following attractive collection of closure properties.
\begin{thmc} For all natural numbers $c,d$ and $n$, the following hold
\begin{itemize}
\item[(1)] $\mathfrak{F}_d$ and $\mathfrak{M}_d$ are subgroup closed
\item[(2)] $\mathfrak{F}_d$ and $\mathfrak{M}_d$  are closed under taking quotients by finite normal subgroups
\item[(3)] $ \mathrm{Ext}(\mathfrak{L}\cap\mathfrak{M}_{c},\mathfrak{M}_{d})\subseteq \mathfrak{M}_{c+d} $
\item[(4)] $ \mathrm{H}_1(n)\mathfrak{F}_d \subseteq  \mathfrak{F}_{n+d}$
\item[(5)] $ \widetilde{\mathrm{H}}_1(n)\mathfrak{F}_d \subseteq  \mathfrak{F}_{n+d}$
\item[(5)] $ \widetilde{\mathrm{H}}_1(n)\mathfrak{M}_d \subseteq  \mathfrak{M}_{n+d}$.
\end{itemize}
\end{thmc}

\section{Preliminaries}
\subsection*{Modules over a category and cohomology}

We start by recalling some of the elementary notions needed to develop a theory of modules over a category and cohomology of a category. For this we will follow the approach of \cite[Ch. 9]{Luck},  to which we refer for further information and proofs of the statements made in this subsection. \\

Let $\mathcal{C}$ be a small category. A right (left) $\mathcal{C}$-module is a contravariant (covariant) functor from $\mathcal{C}$ to the category of abelian groups $\mathrm{Ab}$. The right (left) $\mathcal{C}$-modules are the objects of a category $\mbox{Mod-}\mathcal{C}$ ($\mathcal{C}\mbox{-Mod}$), where the morphism are natural transformations of functors. The category   $\mbox{Mod-}\mathcal{C}$ ($\mathcal{C}\mbox{-Mod}$) is an abelian category, where a sequence
\[0\rightarrow M' \rightarrow M \rightarrow M'' \rightarrow 0\]
is called {\it exact} if it is exact after evaluating in $c$, for all $c \in \mathrm{Ob}(\mathcal{C})$. Take $M \in \mbox{Mod-}\mathcal{C}$ and consider the left exact functor
\[ \mathrm{Hom}_{\mathcal{C}}(M,-) : \mbox{Mod-}\mathcal{C} \rightarrow \mathrm{Ab}: N \mapsto  \mathrm{Hom}_{\mathcal{C}}(M,N), \]
where $\mathrm{Hom}_{\mathcal{C}}(M,N)$ is the abelian group of all natural transformations from $M$ to $N$. By definition,
$M$ is a \emph{projective} $\mathcal{C}$-module if and only if this functor is exact. By considering the contravariant functor  $\mathrm{Hom}_{\mathcal{C}}(-,N)$, one can define injective $\mathcal{C}$-modules in a similar way. Given $M \in \mbox{Mod-}\mathcal{C}$ and $N \in \mathcal{C}\mbox{-Mod}$, one defines their tensor product $M \otimes_{\mathcal{C}}N$ to be the following abelian group
\[ M \otimes_{\mathcal{C}}N = \bigoplus_{c \in \mathrm{Ob}(\mathcal{C})}M(c)\otimes N(c) \big/ I, \]
where $I$ is the abelian group generated by the elements \[M(\varphi)(m)\otimes n - m\otimes N(\varphi)(n),\] for all $\varphi \in \mathrm{Mor}(c,c')$, $m \in M(c')$, $n \in N(c)$ and for all $c,c' \in \mathrm{Ob}(\mathcal{C})$. \\

Given  $M \in \mbox{Mod-}\mathcal{C}$ and $N \in \mathcal{C}\mbox{-Mod}$, one can define the right exact functors
\[ -\otimes_{\mathcal{C}}N : \mbox{Mod-}\mathcal{C}   \rightarrow \mathrm{Ab}: V \rightarrow V \otimes_{\mathcal{C}}N\]
\[ M\otimes_{\mathcal{C}}- :   \mathcal{C}\mbox{-Mod} \rightarrow \mathrm{Ab}: R \rightarrow M \otimes_{\mathcal{C}}R. \]
Next, we describe the free modules in $ \mbox{Mod-}\mathcal{C}$ and $\mathcal{C}\mbox{-Mod}$. \\

Take $c \in  \mathrm{Ob}(\mathcal{C})$ and define the right $\mathcal{C}$-module
\[ \mathbb{Z}[-,c]: \mathcal{C} \rightarrow \mathrm{Ab} : b \rightarrow \mathbb{Z}[b,c], \]
where $\mathbb{Z}[b,c]$ is the free abelian group generated by $\mathrm{Mor}(b,c)$ and a morphism $\varphi \in \mathrm{Mor}(b_1,b_2)$, is mapped to  $$\mathbb{Z}[\varphi,c]: \mathbb{Z}[b_2,c] \rightarrow \mathbb{Z}[b_1,c]$$ where  $\mathbb{Z}[\varphi,c](\alpha)= \alpha \circ \varphi$ for all $\alpha \in \mathrm{Mor}(b_2,c)$. Similarly, one can define the left $\mathcal{C}$-module $\mathbb{Z}[c,-]$ for every  $c \in \mathrm{Ob}(c)$. Arbitrary direct sums of modules of the form $\mathbb{Z}[-,c]$ ($\mathbb{Z}[c,-]$) are called \emph{free} right (left) $\mathcal{C}$-modules. The summands are called \emph{basic free modules}. A free module is called \emph{finitely generated} if it is a direct of finitely many basic free modules.
Free modules are free in the following sense: given a right $\mathcal{C}$-module, an object $c \in \mathrm{Ob}(\mathcal{C})$ and an element $m \in M(c)$, there exists a unique natural transformation $\varphi: \mathbb{Z}[-,c] \rightarrow M $ such that $\varphi(c)(\mathrm{Id_c})=m$, and similarly for left modules. This implies that any $\mathcal{C}$-module is a quotient of a free $\mathcal{C}$-module. A $\mathcal{C}$-module is called \emph{finitely generated} if it is quotient of a finitely generated free module. Basic free modules satisfy the following version of Yoneda's lemma. Let $M \in  \mbox{Mod-}\mathcal{C}$ and $N \in \mathcal{C}\mbox{-Mod}$. For each $c \in \mathrm{Ob}(\mathcal{C})$, there are natural isomorphisms
\[ \mathbb{Z}[-,c]\otimes_{\mathcal{C}}N \xrightarrow{\cong} N(c): \varphi\otimes n \mapsto N(\varphi)(n),   \]
\[ M\otimes_{\mathcal{C}}\mathbb{Z}[c,-] \xrightarrow{\cong} M(c): m\otimes \varphi \mapsto M(\varphi)(m)   \]
and
\[ \mathrm{Hom}_{\mathcal{C}}(\mathbb{Z}[-,c],M) \xrightarrow{\cong} M(c): f \mapsto f_c(\mathrm{Id_c}). \]
Using this last isomorphism, one can easily show that free $\orb$-modules are projective. One can also verify that a right $\mathcal{C}$-module is projective if and only if it is a direct summand of a free $\mathcal{C}$-module.\\

Let $\mathcal{D}$ be a small category. Associated to a functor $\pi: \mathcal{D} \rightarrow \mathcal{C}$, one has a \emph{restriction functor}
\[ \mathrm{res}_{\pi}:  \mbox{Mod-}\mathcal{C} \rightarrow  \mbox{Mod-}\mathcal{D} : N \mapsto N \circ \pi, \]
an \emph{induction functor}
\[ \mathrm{ind}_{\pi}:  \mbox{Mod-}\mathcal{D }\rightarrow  \mbox{Mod-}\mathcal{C}: M \mapsto M(?) \otimes_{\mathcal{D}} \mathbb{Z}[-,\pi(?)] \]
and a \emph{coinduction functor}
\[  \mathrm{coind}_{\pi}:  \mbox{Mod-}\mathcal{D} \rightarrow  \mbox{Mod-}\mathcal{C}: M \mapsto \mathrm{Hom}_{\mathcal{D}}\Big(\mathbb{Z}[\pi(?),-],M(?)\Big).                  \]
Of course, one can define these restriction and (co)induction functors for left modules as well. 
The functor $\mathrm{ind}_{\pi}$ is left adjoint to the functor $\mathrm{res}_{\pi}$ and the functor $\mathrm{coind}_{\pi}$ is right adjoint to the functor $\mathrm{res}_{\pi}$. In particular, there are isomorphisms
\begin{equation} \label{eq: adjoint} \mathrm{Hom}_{\mathcal{C}}(\mathrm{ind}_{\pi}(M),N)\cong \mathrm{Hom}_{\mathcal{D}}(M,\mathrm{res}_{\pi}(N)) \end{equation}
\begin{equation} \label{eq: adjoint2}  \mathrm{Hom}_{\mathcal{D}}(\mathrm{res}_{\pi}(N),M)\cong \mathrm{Hom}_{\mathcal{C}}(N,\mathrm{coind}_{\pi}(M)) \end{equation}
that are natural in $M \in  \mbox{Mod-}\mathcal{D}$ and $N \in  \mbox{Mod-}\mathcal{C}$. In a similar spirit, one also has the following natural isomorphisms
\begin{equation} \label{eq: tensor id1}
\mathrm{ind}_{\pi}(M) \otimes_{\mathcal{C}} N \cong M \otimes_{\mathcal{D}} \mathrm{res}_{\pi}(N)
\end{equation}
\begin{equation} \label{eq: tensor id2}
\mathrm{res}_{\pi}(N) \otimes_{\mathcal{D}} M \cong N \otimes_{\mathcal{C}} \mathrm{ind}_{\pi}(M).
\end{equation}

For each $M \in \mbox{Mod-}\mathcal{C} $, there exists a free right $\mathcal{C}$-module $F$, together with a surjection $F \rightarrow M$. By using the splicing technique, this shows that every $M \in  \mbox{Mod-}\mathcal{C} $ admits a free (and hence projective) resolution. One can also show that every module admits an injective coresolution. Hence, all the ingredients are in place to construct Ext-functors $\mathrm{Ext}^{n}_{\mathcal{C}}(-,-)$ that have all the usual properties.\\

Finally, let us mention that if $\mathcal{C}$ is a preadditive category for which the morphism sets are free abelian groups, then one can develop an entirely similar theory as above by considering the categories of contravariant and covariant additive functors from $\mathcal{C}$ to $\mathrm{Ab}$. Essentially the only difference is that in order to construct basis free modules, one does not have to consider free abelian groups generated by morphism sets, but can use the morphism sets themselves (which already are free abelian). This situation will come up when we consider the Mackey category.

\subsection*{Bredon cohomology}
We will now apply the general theory presented above to the orbit category with finite isotropy. For more details on Bredon (co)homology one can for example consult \cite{FluchThesis} and  \cite{nucinkis04}.
\\

A family of subgroups $\mathcal{F}$ of a group is a collection of subgroups of that group  which is closed under conjugation and taking subgroups. In what follows, $\mathcal{F}$ will always denote the family of finite subgroups of a given group. Let $G$ be a group. The \emph{orbit category} $\orb$ is a category defined as follows: the objects are the transitive $G$-sets $G/H$ for all $H \in \mathcal{F}$ and the morphisms are all $G$-equivariant maps between the objects. Note that every morphism $\varphi: G/H \rightarrow G/K$ is completely determined by $\varphi(H)$, since $\varphi(xH)=x\varphi(H)$ for all $x \in G$. Moreover, there exists a morphism $$G/H \rightarrow G/K : H \mapsto xK \mbox{ if and only if } x^{\scriptscriptstyle -1}Hx \subseteq K.$$ We denote the morphism $\varphi: G/H \rightarrow G/K: H\mapsto xK$  by $G/H \xrightarrow{x} G/K$.

When applying the theory of the previous section to the orbit category $\orb$, we will abbreviate $\mathrm{Hom}_{\orb}(M,N)$ by  $\mathrm{Hom}_{\mathcal{F}}(M,N)$ and $M\otimes_{\orb}N$ by $M\otimes_{\mathcal{F}}N$. The \emph{$n$-th Bredon cohomology of $G$} with coefficients $M \in \orbmod$ is by definition
\[ \mathrm{H}^n_{ \mathcal{F}}(G,M)= \mathrm{Ext}^{n}_{\orb}(\underline{\mathbb{Z}},M) \]
where the functor $\underline{\mathbb{Z}}$ is constantly equal to $\mathbb{Z}$ on objects and maps all morphisms to the identity map. The augmented cellular chain complex $C_{\ast}(\underline{E}G^{-})\rightarrow \underline{\mathbb{Z}}$ of a classifying space for proper actions of $G$ is a free $\orb$-resolution of $\underline{\mathbb{Z}}$ that can be used to compute $ \mathrm{H}^n_{ \mathcal{F}}(G,M)$. The \emph{Bredon cohomological dimension} of $G$ for the family $\mathcal{F}$ is defined as
\[ \underline{\mathrm{cd}}(G) = \sup\{ n \in \mathbb{N} \ | \ \exists M \in \orbmod :  \mathrm{H}^n_{\mathcal{F}}(G,M)\neq 0 \}. \]
The invariant $\underline{\mathrm{cd}}(G)$ is also called the \emph{Bredon cohomological dimension for proper actions}.\\

When $N$ is a subgroup of $G$, we have an obvious functor 
\[ \iota_N^G:  \mathcal{O}_{\mathcal{F}}N \rightarrow \orb : N/H \mapsto G/H.  \]
 The restriction and (co)induction functors associated to this functor will be denoted by $\mathrm{res}_N^G$ and $\mathrm{(co)ind}_N^G$. Note that in this context, we will usually omit restriction functors from the notation. It is well-known that (e.g. see \cite[Lemma 2.9]{Symonds05}) the functor $\mathrm{ind}_N^G$ is exact and that the functor $\mathrm{res}_N^G$ preserves free and hence projective objects. Using this, one obtains a version of Shapiro's lemma for Bredon cohomology:  for each $M \in \mbox{Mod-}\mathcal{O}_{\mathcal{F}}N$ we have
\[ \mathrm{H}^{n}_{\mathcal{F}}(N,M)\cong  \mathrm{H}^{n}_{\mathcal{F}}(G,\mathrm{coind}_N^G(M)) \]
for every $n \in \mathbb{N}$. It follows from this isomorphism that $\underline{\mathrm{cd}}(N)\leq \underline{\mathrm{cd}}(G)$. \\

It is clear that one can also define an orbit category, and therefore Bredon cohomology, for any other family of subgroups of $G$.
\subsection*{Mackey functors}

We now turn to Mackey functors. Mackey functors were introduced for finite groups by Dress and Green, and were studied extensively in this context by Th{\'e}venaz, Webb and others (e.g.~see \cite{ThevenazWebb95}, \cite{Webb00}). However, many of the elementary results obtained about Mackey functors for finite groups generalize to infinite groups and their family of finite subgroups. Our treatment of Mackey functors follows the approach of Mart{\'\i}nez-P{\'e}rez and Nucinkis in \cite{MartinezNucinkis06}. \\

Let $G$ be a group and let $\mathcal{F}$ be the family of finite subgroups of $G$. Consider the diagrams of the form
\begin{equation} \label{eq: mackey morphism}G/S \xleftarrow{\varphi} \Delta \xrightarrow{\psi} G/K \end{equation}
where the maps $\varphi$ and $\psi$ are $G$-equivariant, $S,K \in \mathcal{F}$ and $\Delta = \coprod_{i=1}^n G/H_i$ is a finitely generated $G$-set with stabilizers in $\mathcal{F}$. A diagram of the form $(\ref{eq: mackey morphism})$ is called \emph{equivalent} to a diagram \[G/S \xleftarrow{\varphi'} \Delta' \xrightarrow{\psi'} G/K\] if there exists a $G$-equivariant bijection $\theta: \Delta \rightarrow \Delta'$ such that
\[ \xymatrix{    & \ar[dl]_{\varphi} \Delta \ar[dd]^{\theta} \ar[dr]^{\psi} &  \\
G/S & & G/K \\
& \ar[ul]^{\varphi'} \Delta'  \ar[ur]_{\psi'} &} \]
commutes.
By $\omega_{\mathcal{F}}(S,K)$ we mean the set of equivalence classes of diagrams of the form $(\ref{eq: mackey morphism})$. Elements of $\omega_{\mathcal{F}}(S,K)$ are called \emph{morphisms}. Note that we also allow the \emph{empty} morphism $G/S \leftarrow \emptyset \rightarrow G/K$. One can check that $\omega_{\mathcal{F}}(S,K)$ is a free abelian monoid with disjoint union of $G$-sets as addition and the empty morphism as neutral element. A \emph{basic morphism} in $\omega_G(S,K)$ is a morphism that can be represented by a diagram of the from $G/S \xleftarrow{x} G/L \xrightarrow{y} G/K$. The \emph{Mackey Category} (also called \emph{Burnside category}) $\mathcal{M_F}G$ is defined as follows. Its objects are the transitive $G$-sets $G/H$ for all $H \in \mathcal{F}$. The space of morphisms $\mathrm{Mor}(G/S,G/K)$ is by definition the abelian group completion of $\omega_{\mathcal{F}}(S,K)$. Composition is defined by pullback in the category of $G$-sets on basis morphisms and then extended via linearity to all morphisms. To be more specific: let $\varphi \in \mathrm{Mor}(G/S,G/K)$ be represented by the diagram $G/S \xleftarrow{x} G/L \xrightarrow{x'} G/K$ and let $\psi \in \mathrm{Mor}(G/K,G/H)$ be represented by the diagram $G/K \xleftarrow{y} G/P \xrightarrow{y'} G/H$. Then the composition $\psi \circ \varphi \in \mathrm{Mor}(G/S,G/H)$ is represented by the outer diagram below, obtained by taking the pullback in the category of $G$-sets of $G/L \xrightarrow{x'} G/K \xleftarrow{y} G/P$.
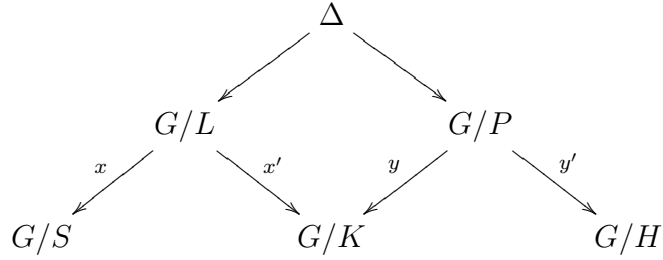
\[ \xymatrix{
&                          & \ar[dl] \Delta  \ar[dr]& & \\
& \ar[dl]_{x} G/L  \ar[dr]^{x'} & & \ar[dl]_{y} G/P \ar[dr]^{y'} &   \\
G/S & & G/K &  & G/H} \]
We will denote the set of morphism from $G/S$ to $G/K$ in $\mathcal{M_F}G$ by $\mathbb{Z}^G[S,K]$. Let us point out that this group is different from $\mathbb{Z}[G/S,G/K]$, which is the free abelian group generated by the morphism from $G/S$ to $G/K$ in the orbit category $\orb$. The category $\mathrm{Mack}_{\mathcal{F}}G$ is the category with objects the contravariant additive functors $M: \mathcal{M_F}G \rightarrow \mathrm{Ab}$, and morphisms all natural transformations between these functors. An object of $\mathrm{Mack}_{\mathcal{F}}G$ is called a \emph{Mackey functor}. Of course one can also consider the category of covariant additive functors from $\mathcal{M_F}G $  to $\mathrm{Ab}$, which is isomorphic to $\mathrm{Mack}_{\mathcal{F}}G$. \\

When applying the theory of the previous section to the Mackey category $\mathcal{M_F}G$ (which is a preadditive category), we will abbreviate $\mathrm{Hom}_{\mathcal{M_F}G}(M,N)$ by  $\mathrm{Hom}_{\mathcal{M_F}}(M,N)$ and $M\otimes_{\mathcal{M_F}G}N$ by $M\otimes_{\mathcal{M_F}}N$. One can construct a functor $\pi$ from the orbit category $\orb$ to the Mackey category $\mathcal{M_F}G$. Obviously, $\pi$ will map an object $G/K$ in $\orb$ to the object $G/K$ in $\mathcal{M_F}G$. A morphism $G/S \xrightarrow{x} G/K$ is by definition mapped to the morphism in the Mackey category represented by $G/S \leftarrow G/S \xrightarrow{x} G/K$. The associated restriction and induction functors 
 \[ \mathrm{res}_{\pi}: \mathrm{Mack}_{\mathcal{F}}G \rightarrow \orbmod: M \mapsto M\circ \pi = M^{\ast}, \] and
\[ \mathrm{ind}_{\pi}: \orb \rightarrow \mathrm{Mack}_{\mathcal{F}}G : N \mapsto N(?) \otimes_{\mathcal{F}}\mathbb{Z}^G[-,\pi(?)] \]
are connected via an adjointness isomorphism
\begin{equation}\label{eq: important adjointness} \mathrm{Hom}_{\mathcal{F}}(N,M^{\ast}) \cong \mathrm{Hom}_{\mathcal{M_F}}(\mathrm{ind}_{\pi}(N),M). \end{equation}
One can show that $ \mathrm{ind}_{\pi}(\mathbb{Z}[-,G/H])=\mathbb{Z}^G[-,H]$ (see \cite[Th. 3.7]{MartinezNucinkis06}). \\

\begin{remark}\label{remark: subclass} \rm The functor \[\mathrm{res}_{\pi}: \mathrm{Mack}_{\mathcal{F}}G \rightarrow \orbmod: M \mapsto  M^{\ast}\]is not a faithfull functor, meaning that there exist Mackey funtors $M_1$ and $M_2$ such that $M_1 \neq M_2$ but $M_1^{\ast}= M_1^{\ast}$. Therefore, the notation $\mathrm{Mack}_{\mathcal{F}}(G) \subseteq \mbox{Mod-}\mathcal{O}_{\mathcal{F}}G$ from the introduction is a little misleading, and we were careful not to say that $\mathrm{Mack}_{\mathcal{F}}(G)$ is a subcategory of $\mbox{Mod-}\mathcal{O}_{\mathcal{F}}G$, but used the `meaningless' term subclass. The notation $\mathrm{Mack}_{\mathcal{F}}(G) \subseteq \mbox{Mod-}\mathcal{O}_{\mathcal{F}}G$ is just meant to point out that one can consider Bredon cohomology with Mackey functor coefficients via the assignment $M \mapsto M^{\ast}$. One could also define cohomology of Mackey functors, within the  category $\mathrm{Mack}_{\mathcal{F}}G$, as the Ext-functors $\mathrm{Ext}^{n}_{\mathcal{M}_{\mathcal{F}}G}(\mathbb{Z}^{G},M)$ where $\mathbb{Z}^{G}=\mathbb{Z}^G[-,G]=\mathrm{ind}_{\pi}(\underline{\mathbb{Z}})$ is the Burnside ring functor. However, Mart{\'\i}nez-P{\'e}rez and Nucinkis prove in \cite[Th. 3.8.]{MartinezNucinkis06} that \[\mathrm{Ext}^{n}_{\mathcal{M}_{\mathcal{F}}G}(\mathbb{Z}^{G},M)\cong \mathrm{Ext}^{n}_{\mathcal{O}_{\mathcal{F}}G}(\underline{\mathbb{Z}},M^{\ast}) ,\] so this does not give anything new.

\end{remark}

In what follows we will also require the following notation. Let $M \in \mathrm{Mack}_{\mathcal{F}}G$ be a Mackey functor. For $H \subseteq K \in \mathcal{F}$ and $g \in G$ one defines the maps
\begin{itemize}
\item[-]  $R_H^K= M(G/H \leftarrow G/H \rightarrow G/K): M(G/K) \rightarrow M(G/H)$
\item[-]  $I_H^K= M(G/K \leftarrow G/H \rightarrow G/H): M(G/H) \rightarrow M(G/K)$
\item[-] $c_g(H)= M(G/{}^gH \xleftarrow{g^{-1}} G/H \rightarrow G/H): M(G/H) \rightarrow M(G/{}^{g}H)$.
\end{itemize}
$R_H^K$ is called a \emph{restriction map}, $I_H^K$ is called an \emph{induction map} or \emph{transfer map} and $c_g(H)$ is called a \emph{conjugation map}. One could also define a Mackey functor as collection of abelian group $\{M(G/H)\}_{H \in \mathcal{F}}$, together with a collection of restriction, induction  and conjugation maps between these abelian groups, that satisfy certain compatibility conditions including the so-called \emph{Mackey axiom} (see Section 2 of \cite{ThevenazWebb95} and \cite{Webb00}). \\

Interesting examples of Mackey functors arise naturally from representation rings, Burnside rings, topological K-theory of classifying spaces and algebraic K-theory of group rings of finite groups, to name just a few. A lot of interesting Mackey functors have an additional property called \emph{conjugation invariance} (e.g. see \cite{HambletonTaylorWillimas10}), which says the the conjugation maps $c_g:  M(G/H) \rightarrow M(G/H)$ are identity maps for all $H \in \mathcal{F}$ and all $g$ contained in the centralizer $C_G(H)$ of $H$ in $G$. For this reason, the conjugation invariance condition is sometimes (e.g. in \cite[Sec. 5]{Luck5}) added to the definition of a Mackey functor. We do not require this, and work with the most general definition of a Mackey functor.\\

A Mackey functor $M \in \mathrm{Mack}_{\mathcal{F}}G$ is called \emph{cohomological} if \[I_H^K \circ R_H^K = [K:H]\mathrm{Id}\] for all $H \subseteq K \in \mathcal{F}$. The subcategory of  $\mathrm{Mack}_{\mathcal{F}}G$  containing all cohomological Mackey functors will be denoted by  $\mathrm{coMack}_{\mathcal{F}}G$.  Note that quotient and sub-Mackey functors of cohomological Mackey functors are again cohomological Mackey functors. Important examples of cohomological Mackey functors include fixed point functors and coinvariance functors. Let $M$ be a left $G$-module. The \emph{fixed point functor} $\underline{M}$ is defined by
the collection of abelian groups $\underline{M}(G/H)=M^{H}$ for all $H \in \mathcal{F}$ where \[M^H=\{ m \in M \ | \ h\cdot m =m \ \mbox{for all}\ h \in H \},\] and for all $L \subseteq K \in \mathcal{F}$ and $g \in G$, the restriction, induction and conjugation maps 
\begin{eqnarray*} R_L^K : M^K &\rightarrow & M^L : m \mapsto m, \\
I_L^K : M^L & \rightarrow & M^K : m \mapsto \sum_{k \in K/L} k\cdot m \end{eqnarray*}
and
\[   c_g(K): M^K \rightarrow M^{{}^gK}: m \mapsto g\cdot m.   \]
It is an easy matter to verify that this data indeed produces a cohomological Mackey functor $\underline{M}: \mathcal{M}_{\mathcal{F}}G \rightarrow \mathrm{Ab}: G/H \mapsto M^H$. The subcategory of $\mathrm{coMack}_{\mathcal{F}}G$ consisting of fixed point functors will be denoted by $\mathcal{F}\mathrm{ix}(G)$. The functor from the category of left $G$-modules to the category of Mackey functors of $G$ that assigns to a $G$-module $M$ the fixed point functor $\underline{M}$ is right adjoint (see \cite[Lemma 4.2]{MartinezNucinkis06}) to the \emph{evaluation at the identity} functor \[\mathrm{ev}: \mathrm{Mack}_{\mathcal{F}}G \rightarrow G\mbox{-Mod}: M \mapsto M(G/e).\] In particular, for every $G$-module $M$ and every Mackey functor $N$ of $G$, there is a natural isomorphism
\[  \mathrm{Hom}_G(N(G/e),M) \xrightarrow{\cong} \mathrm{Hom}_{\mathcal{M}_{\mathcal{F}}}(N,\underline{M}).         \]
Since for any right $\orb$-module $L$, one has $\mathrm{ind}_{\pi}(L)(G/e)\cong L(G/e)$  as $G$-modules (because $\mathbb{Z}^G[e,H]=\mathbb{Z}[G/e,G/H]$), one can combine the isomorphism above with the isomorphism (\ref{eq: important adjointness}) to obtain an adjointness isomorphism
\[  \mathrm{Hom}_G(L(G/e),M) \xrightarrow{\cong} \mathrm{Hom}_{\mathcal{F}}(L,\underline{M}^{\ast})        \]
for every left $G$-module $M$ and every right $\orb$-module $L$. 
One can view $\underline{M}$ as being the zeroth cohomology functor $G/H \mapsto \mathrm{H}^0(H,M)$. More generally, it  turns out that all cohomology functors $G/H \mapsto \mathrm{H}^n(H,M)$ have the structure of a cohomological Mackey functor, as do the homology functors $G/H \mapsto \mathrm{H}_n(H,M)$ (e.g. see \cite[Prop. 9.5]{Brown}). We will call the zeroth homology functor $G/H \mapsto \mathrm{H}_0(H,M)$, the \emph{coinvariance functor} associated to $M$ and denote it by $\overline{M}$.  \\

Suppose that $N$ is a subgroup of $G$. Then, we have a functor $\iota_N^G : \mathcal{M}_{\mathcal{F}\cap N}N \rightarrow \mathcal{M}_{\mathcal{F}}G$ defined in the obvious way. Just as before but now for Mackey functors, this functor induces a restriction functor $ \mathrm{res}_N^G$ that preserves projective objects (see \cite[Prop. 3.1]{MartinezNucinkis06}) and an exact (co)induction functor $\mathrm{(co)ind}_N^G$ satisfying the usual adjointness properties. The restriction functors of the form $\mathrm{res}_N^G$ will usually be omitted from the notation. The following two propositions contain key properties of Mackey functors, which fail in general for right $\mathcal{O}_{\mathcal{F}}G$-modules. These properties are the main reason why Mackey functors are interesting from a cohomological point and view and they explain why Mackeys functor are perhaps a more natural generalization of ordinary $G$-modules than right $\orb$-modules. 
\begin{proposition}{\rm \cite[Lemma 3.2 and Th. 3.3]{MartinezNucinkis06}}\label{prop: key mackey1} Let $N$ be subgroup of $G$ and $M \in \mathrm{Mack}_{\mathcal{F}}N$. 
There is a natural inclusion of Mackey functors \[i:\mathrm{ind}_N^G(M) \rightarrow   \mathrm{coind}_N^G(M).\] Moreover, for each $S \in \mathcal{F}$, we have
\[    \mathrm{ind}_N^G(M)(G/S) \cong \bigoplus_{x \in [N \setminus G/S]}M(N/{}^xS\cap N)   \] \[  \mathrm{coind}_N^G(M)(G/S) \cong \prod_{x \in [N \setminus G/S]}M(N/{}^xS\cap N)     \]
such that the inclusion $i(G/S)$ corresponds to the natural inclusion of the sum into the product. If these natural inclusions are isomorphisms for each $S \in \mathcal{F}$, then $\mathrm{ind}_H^G(M)$ and  $\mathrm{coind}_H^G(M)$ are isomorphic as Mackey functors. In particular, this happens when $[G:N]<\infty$.
\end{proposition}
Using the adjointness isomorphism (\ref{eq: important adjointness}) one can check that for every $M \in \mathrm{Mack}_{\mathcal{F}}N$, there is an isomorphism of right $\orb$-modules \begin{equation}\label{eq: Mackey-bredon coind}\mathrm{coind}_N^G(M)^{\ast} \cong \mathrm{coind}_N^G(M^{\ast}).\end{equation} On the other hand, the right $\orb$-modules $\mathrm{ind}_N^G(M)^{\ast}$ and  $\mathrm{ind}_N^G(M^{\ast})$ are quite different. The reason for this discrepancy and for the fact the that the proposition above is not true for modules over the orbit category is basically that the space of morphisms between objects in the Mackey category is symmetric (i.e.~$\mathbb{Z}^G[H,K]= \mathbb{Z}^G[K,H]$), while morphism sets in the orbit category are highly asymmetrical. \\

To state the second proposition, we need to generalize our notion of Mackey functor in the following way. Let $\mathcal{S}$ be a family of subgroups of $G$. In (\ref{eq: mackey morphism}),  we allow more generally all diagrams of the form
\[G/S \xleftarrow{\varphi} \Delta \xrightarrow{\psi} G/K \]
where the maps $\varphi$ and $\psi$ are $G$-equivariant, $S,K \in \mathcal{S}$ and $\Delta = \coprod_{i=1}^n G/H_i$ is a finitely generated $G$-set such that $H^{x_i}$, where $\varphi(eH_i)=x_iS$, is a finite index subgroup of $S$ for each $i \in \{1,\ldots,n\}$. Now can proceed similarly as before to obtain a category $\mathcal{M}_{\mathcal{S}}G$ which allows for isotropy in all subgroup of $\mathcal{S}$. 
A Mackey functor of $G$ for the family $\mathcal{S}$ is a contravariant additive functor from $\mathcal{M}_{\mathcal{S}}G$ to $\mathrm{Ab}$. As before one can also define restriction, induction and conjugation maps for these Mackey functors, keeping in mind that the induction maps $I_H^K$ will only exist when $H$ is a finite index subgroup of $K$. 

Let $\mathcal{ALL}$ be the family containing all subgroups of $G$. The following proposition was proven in  \cite[Cor. 4.14 and Prop. 5.4]{HambletonPamukYalcin} in the case where $G$ is a finite group and $\mathcal{F}$ is a family of subgroups of $G$. The proof goes trough in the case of infinite groups.
\begin{proposition} {\rm \cite[Cor. 4.14 and Prop. 5.4]{HambletonPamukYalcin}}\label{prop: key mackey2} Let $G$ be a group and let $M \in \mathrm{Mack}_{\mathcal{F}}G$ be a (cohomological) Mackey functor. Then for each $n$,
\[  \mathcal{M}_{\mathcal{ALL}}G \rightarrow \mathrm{Ab} : G/N \mapsto \mathrm{H}^n_{\mathcal{F}}(N,M^{\ast}) \]
is a (cohomological) Mackey functor. 
\end{proposition}
More generally this proposition is also true when $\mathcal{F}$ is replaced by any other family of subgroups containing only finite subgroups.

\section{An upper bound for the Bredon cohomological dimension for proper actions}\label{sec: upper}
Let $G$ be a group and let $\mathcal{F}$ be the family of finite subgroup of $G$. In \cite{nucinkis99} (see also \cite{nucinkis00} and \cite{Gandini12}), Nucinkis introduces a cohomology theory $\mathcal{F}\mathrm{H}^{\ast}(G,M)$ for groups $G$ and $G$-modules $M$, relative to the family of finite subgroups $\mathcal{F}$. The groups $\mathcal{F}\mathrm{H}^{\ast}(G,M)$ can be computed by taking the cohomology of a cochain complex obtained from applying $\mathrm{Hom}_G(-,M)$ to an $\mathcal{F}$-projective resolution $P_{\ast} \rightarrow \mathbb{Z}$. An $\mathcal{F}$-projective resolution $P_{\ast} \rightarrow \mathbb{Z}$ is by definition an exact chain complex of $G$-modules $P_{\ast} \rightarrow \mathbb{Z}$ that is split when restricted to every finite subgroup of $G$, and such that each $P_n$ is a direct summand of a $G$-module of the form $V\otimes \mathbb{Z}[\Delta]$, where $\Delta$ is the $G$-set $\coprod_{H \in \mathcal{F}}G/H$ and $V$ is a $G$-module. The cohomology theory $\mathcal{F}\mathrm{H}^{\ast}(G,-)$ yields the notion of \emph{cohomological dimension of $G$, relative to the family of finite subgroups}
\[  \mathcal{F}\mathrm{cd}(G)=  \sup\{ n \in \mathbb{N} \ | \ \exists M \in G\mbox{-Mod}   :  \mathcal{F}\mathrm{H}^n(G,M)\neq 0 \} .   \]
Using the following known proposition one can also define $ \mathcal{F}\mathrm{cd}(G)$ using Bredon cohomology with fixed point functor coefficients.

\begin{proposition} (\cite[Th. 4.3]{MartinezNucinkis06} or \cite[Eq. (3)]{Martinez07}) \label{prop: rel coh} Let $G$ be a group. Then we have
\[  \mathcal{F}\mathrm{H}^{n}(G,M)\cong \mathrm{H}^{n}_{\mathcal{F}}(G,\underline{M}^{\ast})          \]
for every $G$-module $M$ and every $n \in \mathbb{N}$, and hence
\[  \mathcal{F}\mathrm{cd}(G)=  \sup\{ n \in \mathbb{N} \ | \ \exists \underline{M} \in \mathcal{F}\mathrm{ix}(G)   :  \mathrm{H}^n_{\mathcal{F}}(G,\underline{M}^{\ast})\neq 0 \} .   \]
\end{proposition}
\begin{proof} Let $F$ be a free $\mathcal{O}_{\mathcal{F}}G$-resolution of $\underline{\mathbb{Z}}$. In \cite[Cor. 3.4]{nucinkis00} it is shown that $F(G/e)$ is an $\mathcal{F}$-projective resolution of $\mathbb{Z}$ that can be used to compute  $\mathcal{F}\mathrm{H}^{n}(G,M)$. Since fixed point functors are right adjoint to evaluation at the identity functors, we obtain
 \[\mathrm{H}^{n}_{\mathcal{F}}(G,\underline{M}^{\ast}) =\mathrm{H}^n(\mathrm{Hom}_{\mathcal{F}}(F_{\ast},\underline{M}^{\ast}))\cong  \mathrm{H}^n(\mathrm{Hom}_{G}(F_{\ast}(G/e),M))=  \mathcal{F}\mathrm{H}^{n}(G,M).  \]
\end{proof}
It is an immediate fact that $ \mathcal{F}\mathrm{cd}(G)\leq \underline{\mathrm{cd}}(G)$.

\begin{definition} \rm  Let $H$ be a finite subgroup of a group $G$. An \emph{$l$-chain} of $H$ is a strictly ascending sequence of subgroups
\[    \{e\}=H_0 \subsetneq H_1 \subsetneq \ldots \subsetneq H_{l}=H.      \]
The \emph{length $l(H)$ of a finite subgroup} $H$ of $G$ is the largest integer $l$ such that there exists an $l$-chain of $H$. The \emph{lenght of $G$} is by definition
\[    l(G)=\sup\{ l(H)  \ | \ H \in \mathcal{F}    \}.   \]
Let $M$ be a right $\mathcal{O}_{\mathcal{F}}G$-module. We define the number
\[     \xi(M)=\min\{   l(H) \ | \ H \in \mathcal{F} \ \mbox{such that} \ M(G/H)\neq 0          \}         \]
if $M$ is not the zero module, and put $\xi(0)=\infty$. It is clear that  $\xi(M)\leq l(G)$ for every non-zero right $\mathcal{O}_{\mathcal{F}}G$-module $M$. If the group $G$ has a uniform bound on the orders of its finite subgroups, then $ l(G) < \infty $. On the other hand, there are many groups with $l(G)<\infty$ but unbounded torsion.
\end{definition}
Let $H$ be a finite subgroup of a group $G$. The \emph{normalizer} of $H$ in $G$ is denoted by $N_G(H)$ and the \emph{Weyl-group} $W_G(H)$ is the quotient $N_G(H)/H$.
\begin{lemma} \label{lemma: mackey res} Let $G$ be a group and let $M$ be a non-zero right $\mathcal{O}_{\mathcal{F}}G$-module. For every integer $d\geq 0$ exists an exact sequence of right $\orb$-modules
\[ 0 \rightarrow M \xrightarrow{f_0} M_0^{\ast} \xrightarrow{f_1} M_1^{\ast} \xrightarrow{f_2} \ldots \xrightarrow{f_d} M_d^{\ast} \rightarrow Q \rightarrow 0      \]
such that $d+\xi(M)<\xi(Q)$ and $M_i\in \mbox{Mack}_{\mathcal{F}}(G)$ for all $i \in \{0,\ldots,d\}$. Moreover, one has
\[     \mathrm{H}^{n}_{\mathcal{F}}(G,M_i^{\ast}) \cong   \prod_{\substack{H \in (\mathcal{F}) \ \mbox{\tiny s.t.}\\ l(H)\geq i+\xi(M)}} \mathrm{H}^{n}_{\mathcal{F}}(W_G(H),\underline{K_i(G/H)}^{\ast})\]
for every $i \in \{0,\ldots,d\}$ and each $n \in \mathbb{N}$. Here $K_i$ is the image of $f_i$ and  $(\mathcal{F})$ is a set containing one representative of each conjugacy class in $\mathcal{F}$.
\end{lemma}
\begin{proof} Let $M$ be a non-zero right $\mathcal{O}_{\mathcal{F}}G$-module $M$. Take $H \in \mathcal{F}$ and let $W_G(H)$ be the Weyl-group of $H$ considered as a category with one object. Associated to the functor $W_G(H) \rightarrow \mathcal{O}_{\mathcal{F}}G : W_G(H)\mapsto G/H$, there is a
restriction functor
\[ \mathrm{res}:   \mbox{Mod-}\mathcal{O}_{\mathcal{F}}G \rightarrow  \mbox{Mod-}W_G(H): M \mapsto M(G/H)        \]
and a coinduction functor
\[  \mathrm{coind} : \mbox{Mod-}W_G(H) \rightarrow  \mbox{Mod-}\mathcal{O}_{\mathcal{F}}G : V \mapsto \mathrm{Hom}_{W_G(H)}(\mathbb{Z}[G/H,-],V), \]
such that $\mathrm{coind}$ is the right adjoint of $\mathrm{res}$. Define $D_HM=\mathrm{coind}(M(G/H))$ and note that the adjoint of $\mathrm{Id}: M(G/H) \rightarrow M(G/H)$ yields an
canonical map $j_H:M \rightarrow D_HM$ which is the identity map when evaluated at $H$, and an isomorphism when evaluated at any conjugate subgroup $H^g$. Denote by $(\mathcal{F})$ a subset of $\mathcal{F}$ containing exactly one representative from each conjugacy class in $\mathcal{F}$. Now define $DM=\prod_{H \in (\mathcal{F}) } D_HM$ and notice that the maps $j_H$ assemble to form an inclusion $f_0: M \rightarrow DM$. We construct an exact sequence\[ 0 \rightarrow M \xrightarrow{f_0} DM \rightarrow CM \rightarrow 0 \]
and claim that $f_0(K)$ is an isomorphism for all $K \in \mathcal{F}$ such that $l(K) \leq \xi(M)$. Indeed, take $K \in \mathcal{F}$ such that $l(K)\leq \xi(M)$. If $D_HM(G/K)\neq 0$, then this implies that $H$ is subconjugate to $K$. Assume that $H$ is conjugate to a proper subgroup of $K$. Then we have $l(H)<l(K)$ and therefore $M(G/H)=0$. This however implies that $D_HM(G/K)=0$, which is a contradiction. We conclude that $H$ is conjugate to $K$. In this case $j_H(G/K):M(G/K)\rightarrow D_HM(G/K)$ is an isomorphism. Since $(\mathcal{F})$ contains exactly one $H$ that is conjugate to $K$, the claim is proven.

We deduce that $\xi(CM)$ is strictly larger than $\xi(M)$.
Continuing this process and denoting $C(C^{k-1}M)$ by $C^kM$, we obtain an exact sequence
\[   0 \rightarrow M \xrightarrow{f_0} DM \xrightarrow{f_1} DCM \xrightarrow{f_2}  DC^2M \xrightarrow{f_3} \ldots \xrightarrow{f_d} DC^dM \rightarrow C^{d+1}M \rightarrow 0   \]
such that \[\xi(M)< \xi(CM) < \xi(C^{2}M) < \xi(C^3M) < \ldots < \xi(C^{d+1}M). \] This implies that $\xi(M)+i \leq \xi(C^{i}M)$ for all $i$.
In Example 4.8 of \cite{HambletonPamukYalcin} it is shown that the modules $D(-)$ in the exact sequence
\[   0 \rightarrow M \rightarrow DM \rightarrow DCM \rightarrow  DC^2M \rightarrow \ldots \rightarrow DC^{d}M \rightarrow C^{d+1}M \rightarrow 0     \]
have a Mackey functor structure. Now denote $M_i^{\ast}=DC^{i}M$, $C^{d+1}M=Q$ and $K_i=C^iM$. It remains to show that
\[     \mathrm{H}^{n}_{\mathcal{F}}(G,M_i^{\ast}) \cong   \prod_{\substack{H \in (\mathcal{F}) \ \mbox{\tiny s.t.}\\ l(H)\geq i+\xi(M)}} \mathrm{H}^{n}_{\mathcal{F}}(W_G(H),\underline{K_i(G/H)}^{\ast})\]
for every $i \in \{0,\ldots,d\}$. First note that since $\xi(M)+i \leq \xi(K_i)$ and cohomology commutes with products, it suffices to prove that
\[     \mathrm{H}^{n}_{\mathcal{F}}(G,D_HM) \cong   \mathrm{H}^{n}_{\mathcal{F}}(W_G(H),\underline{M(G/H)}^{\ast})\]
for every right $\mathcal{O}_{\mathcal{F}}G$-module $M$ and any $H \in \mathcal{F}$. Recall that $D_HM=\mathrm{coind}(M(G/H))$ and that $\mathrm{coind}$ is the right adjoint of $\mathrm{res}$. 
Now let $F_{\ast}\rightarrow \underline{\mathbb{Z}}$ be a free $\mathcal{O}_{\mathcal{F}}G$-resolution of $\underline{\mathbb{Z}}$. Then we have,
\begin{eqnarray}
 \mathrm{H}^{n}_{\mathcal{F}}(G,D_HM) & = & \mathrm{H}^n(\mathrm{Hom}_{\mathcal{O}_{\mathcal{F}}G}(F_{\ast},D_HM) )\nonumber\\
& \cong  & \mathrm{H}^{n}(\mathrm{Hom}_{W_G(H)}(F_{\ast}(G/H),M(G/H)) ).\label{eq: first one1}
\end{eqnarray}
Consider the quotient map $\pi: N_G(H) \rightarrow W_G(H)$ and the functor
\[     p :   \mathcal{O}_{\mathcal{F}}W_G(H) \rightarrow  \mathcal{O}_{\mathcal{F}}N_G(H): W_G(H)/L\mapsto N_G(H)/\pi^{-1}(L), \]
and let $\mathrm{res}_p$ be the associated restriction functor. The functor $\mathrm{res}_p$ is exact. Moreover, one easily verifies that $\mathrm{res}_p(\mathbb{Z}[-,N_G(H)/K])=\mathbb{Z}[-,W_G(H)/\pi(K)]$ for every finite subgroup $K$ of $N_G(H)$ that contains $H$ and $\mathrm{res}_p(\mathbb{Z}[-,N_G(H)/K])=0$ for all finite subgroups  $K$ of $N_G(H)$ that do not contain $H$.  Therefore, the functor  $\mathrm{res}_p$ also preserves free and hence projective modules. 
Now view, via restriction, the resolution $F_{\ast}\rightarrow \underline{\mathbb{Z}}$ as a free $\mathcal{O}_{\mathcal{F}}N_G(H)$-resolution of $\underline{\mathbb{Z}}$. Then $\mathrm{res}_p(F_{\ast})\rightarrow \underline{\mathbb{Z}}$ is a free $\mathcal{O}_{\mathcal{F}}W_G(H)$-resolution of $\underline{\mathbb{Z}}$. Therefore, we have
\begin{eqnarray}
 \mathrm{H}^{n}_{\mathcal{F}}(W_G(H),\underline{M(G/H)}^{\ast}) & = &  \mathrm{H}^n(\mathrm{Hom}_{\mathcal{O}_{\mathcal{F}}W_G(H)}(\mathrm{res}_p(F_{\ast}),\underline{M(G/H)}^{\ast})) \nonumber\\
 & \cong &   \mathrm{H}^n(\mathrm{Hom}_{W_G(H)}(\mathrm{res}_p(F_{\ast})(W_G(H)/e),M(G/H)) )\nonumber \\
 & \cong &   \mathrm{H}^n(\mathrm{Hom}_{W_G(H)}(F_{\ast}(G/H),M(G/H))) \label{eq: second one1},
\end{eqnarray}
where the first isomorphism follows from the fact that fixed point functors are right adjoint to evaluation at the identity functor. Combining (\ref{eq: first one1}) and (\ref{eq: second one1}), we conclude that 
\[     \mathrm{H}^{\ast}_{\mathcal{F}}(G,D_HM) \cong   \mathrm{H}^{\ast}_{\mathcal{F}}(W_G(H),\underline{M(G/H)}).\]
This finishes the proof.
\end{proof}
\begin{remark} \rm The existence of the exact sequence in the previous lemma is proven in \cite[Lemma 4.10]{HambletonPamukYalcin} for finite groups. Our proof is an immediate generalization of that proof. The observation that\[     \mathrm{H}^{n}_{\mathcal{F}}(G,M_i^{\ast}) \cong   \prod_{\substack{H \in (\mathcal{F}) \ \mbox{\tiny s.t.}\\ l(H)\geq i+\xi(M)}} \mathrm{H}^{n}_{\mathcal{F}}(W_G(H),\underline{K_i(G/H)}^{\ast})\] has, as far as we are aware, not appeared in the literature before.
\end{remark}
Using this lemma, we can prove the first part of Theorem A from the introduction.
\begin{theorem}\label{th: length orbit} Let $G$ be a group with $l(G)<\infty$. Then we have
\[    \underline{\mathrm{cd}}(G) \leq \max_{H \in \mathcal{F}}\Big\{  \mathcal{F}\mathrm{cd}(W_G(H))+l(H)  \Big\}.   \]

\end{theorem}
\begin{proof} We may assume that $\mathcal{F}\mathrm{cd}(G) <\infty$, otherwise there is nothing to prove. Let $n=\max_{H \in \mathcal{F}}\Big\{  \mathcal{F}\mathrm{cd}(W_G(H))+l(H)  \Big\}$. To prove the theorem we must show that $\mathrm{H}^{n+1}_{\mathcal{F}}(G,M)=0$ for every right $\mathcal{O}_{\mathcal{F}}G$-module $M$. To this end, fix a right $\mathcal{O}_{\mathcal{F}}G$-module $M$ and assume by contradiction that $\mathrm{H}^{n+1}_{\mathcal{F}}(G,M)\neq 0$ . Now let 
\[ 0 \rightarrow M \xrightarrow{f_0} M_0^{\ast} \xrightarrow{f_1} M_1^{\ast} \xrightarrow{f_2} \ldots \xrightarrow{f_{d}} M_d^{\ast} \rightarrow 0      \]
be the exact sequence constructed in Lemma \ref{lemma: mackey res} for $d=l(G)$ (this forces $Q=0$). Recall that the image of $f_i$ is denoted by $K_i$. By considering the long exact cohomology sequences associated to
\[  0 \rightarrow K_{i} \rightarrow  M^{\ast}_i \xrightarrow{f_{i+1}} K_{i+1} \rightarrow 0 \]
and noticing that \[\mathrm{H}^{n+1-i}_{\mathcal{F}}(G,M^{\ast}_i)=  \prod_{\substack{H \in (\mathcal{F})\ \mbox{\tiny s.t.}\\ l(H)\geq i+\xi(M)}} \mathrm{H}^{n+1-i}_{\mathcal{F}}(W_G(H),\underline{K_i(G/H)}^{\ast}) =0\]  for all $i\leq d$ (since $\mathcal{F}\mathrm{cd}(W_G(H))<n+1-i$, if $l(H)\geq i$), it follows that
\[  \mathrm{H}^{n+1}_{\mathcal{F}}(G,M)=\mathrm{H}^{n+1}_{\mathcal{F}}(G,K_0) \neq 0 \Longrightarrow \mathrm{H}^{n}_{\mathcal{F}}(G,K_1) \neq 0 \Longrightarrow  \mathrm{H}^{n-1}_{\mathcal{F}}(G,K_2) \neq 0  \Longrightarrow    \]
\[ \ldots \Longrightarrow  \mathrm{H}^{n+1-i}_{\mathcal{F}}(G,K_i) \neq 0 \Longrightarrow \ldots \Longrightarrow \mathrm{H}^{n+1-d}_{\mathcal{F}}(G,M^{\ast}_d) =\mathrm{H}^{n+1-d}_{\mathcal{F}}(G,K_d) \neq 0.\]
This is a contradiction. We conclude that $\mathrm{H}^{n+1}_{\mathcal{F}}(G,M)=0$, and hence $ \underline{\mathrm{cd}}(G)\leq n$.
\end{proof}
\begin{remark} \label{remark: attained upper bounds} \rm
By combining Theorem C of \cite{DP2} together with Example 3.6 of \cite{Martinez12}, as Example 6.5 of \cite{DP2}, one obtains for every natural number $s \in \mathbb{N}$ a virtually torsion free group $G_s$ such that $l(G_s)=s$ and 
\[\mathcal{F}\mathrm{cd}(G_s)+s=\underline{\mathrm{cd}}(G_s).\]
This could cause one to be sceptical regarding Nucinkis' conjecture for groups with $l(G)=\infty$. However, one should keep in mind that in the example mentioned above (and all other known examples of a similar nature), $\mathcal{F}\mathrm{cd}(G_s)$ grows linearly with $s$. So it is still conceivable that there exists an upper bound for $\underline{\mathrm{cd}}(G)$ in terms of  $\mathcal{F}\mathrm{cd}(G)$ that does not involve the notion of length of finite subgroups. For example, there are no counterexamples known to the inequality $\underline{\mathrm{cd}}(G)\leq \mathcal{F}\mathrm{cd}(G)^2 $. Let us mention that it is also still an open problem whether $\underline{\mathrm{cd}}(G)= \mathcal{F}\mathrm{cd}(G)$ when $G$ admits a cocompact model for $\underline{E}G$ or an $\orb$-resolution of $\underline{\mathbb{Z}}$ consisting of finitely generated free $\orb$-modules. The groups $G_s$ above do not admit such a model or such a resolution.

\end{remark}

The following corollary is immediate.
\begin{corollary} \label{cor: bounds}Let $G$ be a group of finite length such that 
\[    \mathcal{F}\mathrm{cd}(W_G(H))+l(H) \leq  \mathcal{F}\mathrm{cd}(G)    \]
for every finite subgroup $H$ of $G$. Then we have
\[ \mathcal{F}\mathrm{cd}(G)=\underline{\mathrm{cd}}(G).\]

\end{corollary}
We end this section with the following theorem, which is the most general statement we can prove regarding the Kropholler-Mislin/Nucinkis  conjecture.
\begin{theorem} \label{th: krophollermislin} Let $G$ be a group. The following three statements are equivalent.
\begin{itemize}
\item[(i)]  There exists a finite dimensional contractible proper $G$-CW-complex. There exists an integer $d$ and a finite dimensional proper $G$-CW complex $X$ such that $X^{H}$ is contractible for all finite subgroups $H$ of $G$ for which $l(H)\geq d$.

\item[(ii)] There exist integers $m$ and $d$ such that $\mathcal{F}cd(G) \leq m$ and $\mathrm{H}^{m+1}_{\mathcal{F}}(G,M)=0$ for all $M$ with $\xi(M)\geq d$.

\item[(iii)] There exists a finite dimensional model for $\underline{E}G$.

\end{itemize}

\end{theorem}
\begin{proof} First assume that $(i)$ holds. As mentioned in the introduction, $n_G<\infty$ implies  that $\mathcal{F}cd(G)<\infty$. Let $m$ be the dimension of $X$. Note that we may assume that $\mathcal{F}cd(G) \leq m$ by replacing $X$ with the join $X \ast \Big( \coprod_{H \in \mathcal{F}}G/H\Big)$ a sufficient amount of times. By applying the Equivariant Extension Theorem (\cite[Prop. 2.3]{Luck}) to the $G$-map $X \rightarrow G/G$ and the function 
\[\eta: \mathcal{F} \rightarrow \mathbb{Z} \cup \{\infty\} : H \mapsto \left\{  \begin{array}{cc} \infty & \mbox{if $l(H)\geq d$} \\ -1 & \mbox{otherwise,} \end{array}\right.\]
we obtain a model $Y$ for $\underline{E}G$ such that all cells above dimension $m$ have isotropy groups $H$ with $l(H)< d$. Let $M$ be a right $\mathcal{O}_{\mathcal{F}}G$-module for which $\xi(M)\geq d$. By definition, we have $M(G/H)=0$ for a $H \in \mathcal{F}$ for which $l(H)< d$. Because $Y$ is a model for $\underline{E}G$, it follows that
\[  \mathrm{H}^i_{\mathcal{F}}(G,M) \cong \mathrm{H}^i(\mathrm{Hom}_{\mathcal{F}}(C_{\ast}(Y^{-}),M)).          \]
Since for all $i> m$ one has
\[   \mathrm{Hom}_{\mathcal{F}}(C_{i}(Y^{-}),M)= \mathrm{Hom}_{\mathcal{F}}(\bigoplus_{\alpha \in I}\mathbb{Z}[-,G/H_{\alpha}],M) = \prod_{\alpha \in I}M(G/H_{\alpha})   \]
with $l(H_{\alpha})< d$, we conclude that $\mathrm{H}^{m+1}_{\mathcal{F}}(G,M)=0 $. This shows that $(i)$ implies $(ii)$.  It is trivial that (iii) implies (i), so it remains to show that (ii) implies (iii).

Assume (ii) holds and let $M$ be a right $\orb$-module. By Lemma \ref{lemma: mackey res}, there is an exact sequence of right $\orb$-modules
\[ 0 \rightarrow M \xrightarrow{f_0} M_0^{\ast} \xrightarrow{f_1} M_1^{\ast} \xrightarrow{f_2} \ldots \xrightarrow{f_{d-1}} M_{d-1}^{\ast} \rightarrow Q \rightarrow 0      \]
such that $d \leq \xi(Q)$ and $\mathrm{H}^{r}_{\mathcal{F}}(G,M_i^{\ast}) =0$ for every $i \in \{0,\ldots,d-1\}$ and all $r>m$. Now proceeding as in the proof of Theorem \ref{th: length orbit}, one obtains an implication
\[    \mathrm{H}^{m+d+1}_{\mathcal{F}}(G,M)\neq 0 \Longrightarrow    \mathrm{H}^{m+1}_{\mathcal{F}}(G,Q)\neq 0.    \]
Since $\xi(Q)\geq d$, this proves that $\mathrm{H}^{m+d+1}_{\mathcal{F}}(G,M)=0$ for every right $\orb$-module $M$. We conclude that $\underline{\mathrm{cd}}(G)\leq m+d$, so there exists a finite dimensional model for $\underline{E}G$.
\end{proof}

\section{Cohomological dimensions for Mackey functors}
Let $G$ a group and let $\mathcal{F}$ be the family of finite subgroups of $G$. By only allowing Mackey functors as coefficients for Bredon cohomology, Mart{\'\i}nez-P{\'e}rez and Nuckinkis define the following invariant.
\begin{definition}(see \cite[Def. 3.4]{MartinezNucinkis06} ) \rm  The \emph{Mackey cohomological dimension} of a group $G$ is
\[ \underline{\mathrm{cd}}_{\mathcal{M}}(G) = \sup\{ n \in \mathbb{N} \ | \ \exists M \in \mathrm{Mack}_{\mathcal{F}}(G) :  \mathrm{H}^n_{\mathcal{F}}(G,M^{\ast})\neq 0 \} .                       \]
\end{definition}

It as an immediate consequence of Shapiro's Lemma and (\ref{eq: Mackey-bredon coind}) that \begin{equation}\label{eq: mackey subgroups}\underline{\mathrm{cd}}_{\mathcal{M}}(N) \leq \underline{\mathrm{cd}}_{\mathcal{M}}(G)\end{equation} for every subgroup $N$ of $G$.
By restricting the coefficients even further to cohomological Mackey functors, we obtain yet another notion of cohomological dimension.
\begin{definition} \rm  For a group $G$, we define
\[ \underline{\mathrm{cd}}_{co\mathcal{M}}(G) = \sup\{ n \in \mathbb{N} \ | \ \exists M \in \mathrm{coMack}(G) :  \mathrm{H}^n_{\mathcal{F}}(G,M^{\ast})\neq 0 \} .                       \]
\end{definition}

We now have four notions of cohomological dimension of a group $G$ for the family of finite subgroups, that satisfy
\begin{equation}  \label{eq: five}   \mathcal{F}\mathrm{cd}(G) \leq \underline{\mathrm{cd}}_{co\mathcal{M}}(G) \leq \underline{\mathrm{cd}}_{\mathcal{M}}(G) \leq \underline{\mathrm{cd}}(G).    \end{equation}
In the previous section we have already examined the relationship between $\mathcal{F}\mathrm{cd}(G)$ and $\underline{\mathrm{cd}}(G)$ when $l(G)<\infty$. In \cite[Th. 5.1]{MartinezNucinkis06} Mart{\'i}nez-P{\'e}rez and Nucinkis show that $\mathrm{vcd}(G)=\mathcal{F}\mathrm{cd}(G)=\underline{\mathrm{cd}}_{\mathcal{M}}(G)$ when $G$ is virtually torsion-free. We will show that the equality $\mathcal{F}\mathrm{cd}(G)=\underline{\mathrm{cd}}_{\mathcal{M}}(G)$ holds more generally for groups $G$ of finite length. In the process, we will also show that the invariant $\underline{\mathrm{cd}}_{co\mathcal{M}}(G)$ behaves as expected when passing to subgroups, and is closely related to the invariant $\mathcal{F}\mathrm{cd}(G)$. \\

We start with some technical results.

\begin{lemma} \label{lemma: ind fixed point}Let $G$ be a group and let $N$ be subgroup of $G$. For every $N$-module $M$ one has an isomorphism of Mackey functors

\[\mathrm{coind}_{N}^G(\underline{M}) \cong  \underline{\mathrm{coind}_N^G(M)} \]
and 
\[\mathrm{ind}_{N}^G(\underline{M}) \cong  \underline{\mathrm{ind}_N^G(M)} .\]

\end{lemma}
\begin{proof} Since fixed point functors are right adjoint to evaluation at the identity functors, and coinduction functors are right adjoint to restriction functors, we obtain for every $K \in \mathrm{Mack}_{\mathcal{F}}N$ an isomorphism
\[    \mathrm{Hom}_{\mathcal{M}_{\mathcal{F}}}(K, \mathrm{coind}_{N}^G(\underline{M}) ) \xrightarrow{\cong}   \mathrm{Hom}_{\mathcal{M}_{\mathcal{F}}}(K, \underline{\mathrm{coind}_N^G(M)} )  \]
that is natural in $K$.  By letting $K$ equal the free modules $\mathbb{Z}^G[-,H]$, for all $H \in \mathcal{F}$, and using the naturality in $K$, we obtain the isomorphism
\[\Psi: \mathrm{coind}_{N}^G(\underline{M}) \xrightarrow{\cong}  \underline{\mathrm{coind}_N^G(M)} .\]
Combining this map with the natural inclusion $\mathrm{ind}_{N}^G(\underline{M})\rightarrow \mathrm{coind}_{N}^G(\underline{M})$ from Proposition \ref{prop: key mackey1}, we obtain an inclusion of Mackey functors
\[  \Omega:  \mathrm{ind}_{N}^G(\underline{M}) \rightarrow  \underline{\mathrm{coind}_N^G(M)}    \]
The natural inclusion of $G$-modules $\mathrm{ind}_N^G(M) \rightarrow \mathrm{coind}_N^G(M)$ yields an inclusion of Mackey functors $\Phi: \underline{\mathrm{ind}_N^G(M)}  \rightarrow \underline{\mathrm{coind}_N^G(M)}$.
By considering the evaluation at $G/e$, one sees that the image of $\Omega$ exactly coincides with image of $\Phi$. Hence, there is an isomorphism.
\[    \mathrm{ind}_{N}^G(\underline{M}) \xrightarrow{\cong}  \underline{\mathrm{ind}_N^G(M)}.    \]

\end{proof}

The following proposition is a generalization from finite groups to infinite groups of Theorem 16.5(i) in \cite{ThevenazWebb95}. It gives fixed point functors a special status inside the category of cohomological Mackey functors.
\begin{proposition} \label{prop: fixed point quotient} Every cohomological Mackey functor is a quotient of a fixed point functor.

\end{proposition}

\begin{proof} Let $M \in \mathrm{coMack}_{\mathcal{F}}G$, and let $m \in M(G/H)$ for $H \in \mathcal{F}$. Define the natural transformation
\[ \Theta_{H,m}:  \underline{\mathbb{Z}}   \rightarrow \mathrm{res}_H^G(M)     \]
by setting $\theta(H/L)(1)=R^H_L(m)$. Using the fact that $M$ (and hence  $\mathrm{res}_H^G(M)$) is a cohomological Mackey functor, one checks that $\Theta_{H,m}$ is well-defined. By the adjointness of induction and restriction functors and the second isomorphism of Lemma \ref{lemma: ind fixed point},  we obtain a natural transformation
\[ \Psi_{H,m}: \underline{\mathbb{Z}[G/H]} \rightarrow M \]
such that $\Psi_{H,m}(G/H)(eH)=m$. 

Although this is not necessary for the proof, we will also give an explicit description of $\Psi_{H,m}$. For this we first need to agree on a basis. Let $H,K \in \mathcal{F}$. By viewing $G/H$ as a $K$-set and writing it as a disjoint union of transitive $K$-sets, one checks that there is an isomorphism of $K$-modules
\[     \bigoplus_{x \in [K \setminus G / H] }\mathbb{Z}[K/K\cap {}^xH] \rightarrow \mathbb{Z}[G/H]:  k( K\cap {}^xH) \mapsto kxH.     \]
Since  $\mathbb{Z}[K/K\cap {}^xH]^K$ is an infinite cyclic group generated by $ \sum_{k \in K/ K \cap {}^xH}k(K\cap {}^xH)  $, it follows that 
\[  \Big\{    \sum_{k \in K/ K \cap {}^xH}kx H \ | \ x \in  [K \setminus G / H] \Big\} \]
is a $\mathbb{Z}$-module basis for $ \mathbb{Z}[G/H]^K$. 
Using this basis, the map $\Psi_{H,m}(G/K)$ is given explicitly by
\[    \Psi_{H,m}(G/K):  \mathbb{Z}[G/H]^K \rightarrow M(G/K):   \sum_{k \in K/K\cap {}^xH}kxH \mapsto I_{K\cap {}^xH}^K \circ R_{K \cap {}^xH}^{{}^xH}\circ c_x(H)(m).      \]

Now let $\{m_{i}\}_{i \in I_H}$ be a set of $\mathbb{Z}$-module generators of $M(G/H)$. Then
   \[   \Psi_H=\bigoplus_{i \in I_H}\Psi_{H,m_i} : \bigoplus_{i\in I_H}\underline{\mathbb{Z}[G/H]} \rightarrow M \] 
   is a natural transformation of Mackey functors such that $\Psi_H(G/H)$ is surjective. Hence, the natural transformation
   \[    \Psi= \bigoplus_{H \in \mathcal{F}} \Psi_H :   \bigoplus_{H \in \mathcal{F}}\bigoplus_{i\in I_H}\underline{\mathbb{Z}[G/H]}  \rightarrow M\]
   is surjective.  Since
   \[   \bigoplus_{H \in \mathcal{F}}\bigoplus_{i\in I_H}\underline{\mathbb{Z}[G/H]}  =  \underline{\bigoplus_{H \in \mathcal{F}}\bigoplus_{i\in I_H}\mathbb{Z}[G/H] },\]
   we conclude that $M$ is a quotient of a fixed point functor. 

\end{proof}

To show that the invariant $\underline{\mathrm{cd}}_{co\mathcal{M}}(G)$ behaves well when passing to subgroups, one has to invoke Shapiro's Lemma. This however requires that the coinduction of a cohomological Mackey functor is also a cohomological Mackey functor. Rather than showing this directly, we will prove this by using the proposition above. 

\begin{proposition} If $N$ is a subgroup of $G$, then
\[\underline{\mathrm{cd}}_{co\mathcal{M}}(N) \leq \underline{\mathrm{cd}}_{co\mathcal{M}}(G). \]
\end{proposition}
\begin{proof} Let $R \in \mathrm{Mack}_{\mathcal{F}}N$ be a cohomological Mackey functor. It follows from Proposition \ref{prop: fixed point quotient} that there exists an $N$-module $M$ and a surjection of Mackey functors $\underline{M}\rightarrow R$. By applying $\mathrm{coind}_N^G$ and using the first isomorphism of Lemma \ref{lemma: ind fixed point}, we obtain a surjection of Mackey functors
\[  \underline{\mathrm{coind}_N^G(M)} \rightarrow \mathrm{coind}_N^G(R). \]
Since $\underline{\mathrm{coind}_N^G(M)} $ is a cohomological Mackey functor for $G$, we conclude that $\mathrm{coind}_N^G(R)$ is also a cohomological Mackey functor for $G$. The proposition now follows from Shapiro's lemma and (\ref{eq: Mackey-bredon coind}).
\end{proof}

As announced earlier, we show that the invariant $\underline{\mathrm{cd}}_{co\mathcal{M}}(G)$ is closely related to the invariant $ \mathcal{F}\mathrm{cd}(G)$.

\begin{proposition}  \label{prop: Fcd2}  Let $G$ be a group such that $\underline{\mathrm{cd}}_{co\mathcal{M}}(G)<\infty$, then
\[    \mathcal{F}\mathrm{cd}(G)=\underline{\mathrm{cd}}_{co\mathcal{M}}(G) .  \]

\end{proposition}
\begin{proof} Assume that  $\underline{\mathrm{cd}}_{co\mathcal{M}}(G)=n<\infty$, and let $M$ be a cohomological Mackey functor such that $\mathrm{H}^{n}_{\mathcal{F}}(G,M^{\ast}) \neq 0$. By Proposition \ref{prop: fixed point quotient}, there exists a $G$-module $V$ and a surjection of Mackey functors  $\underline{V} \rightarrow M $. A long exact cohomology sequence then implies that $\mathrm{H}^{n}_{\mathcal{F}}(G,\underline{V}^{\ast}) \neq 0$. We conclude that $ \mathcal{F}\mathrm{cd}(G)=\underline{\mathrm{cd}}_{co\mathcal{M}}(G) $.
\end{proof}

\begin{definition}\rm Let $G$ be a group, let $M \in \mbox{Mack}_{\mathcal{F}}G$ be a Mackey functor and let $F$ be a finite subgroup of $G$. The \emph{subfunctor of $M$ generated by $F$} is the Mackey functor 
\[M_F= \bigcap \Big\{   \mbox{$R$ is a Mackey subfunctor of $M$ for which $R(G/F)=M(G/F)$}   \Big\}. \]
Since the morphism $G/F^{g} \xleftarrow{g} G/F \rightarrow G/F$ in the Mackey category induces an isomorphism between $M(G/F)$ and $M(G/F^{g})$, it follows that $M_F =M_{F^g}$ for every $g \in G$.
One says that $M$ is \emph{generated by a collection $\{F_i\}_{i \in I}$ of finite subgroups} of $G$, if the canonical map 
\[  \bigoplus_{i \in I} M_{F_i} \rightarrow M  \]
is surjective. Note that one may always assume that the collection  $\{F_i\}_{i \in I}$ contains at most one representative of each conjugacy class of finite subgroups of $G$.

\end{definition}

The following lemma is implicitly contained in the proof of Theorem 5.1 in \cite{MartinezNucinkis06}.
\begin{lemma} \label{lemma: inflation}Let $G$ be a group and let $K$ be a finite subgroup of $G$.  Let $M \in \mbox{Mack}_{\mathcal{F}}G$ be a Mackey functor that is generated by $K$ such that
$M(G/S)=0$ for all $S \subsetneq K$. Denote by $Q$ the normalizer of $K$ in $G$. There exists a Mackey functor $T  \in  \mbox{Mack}_{\mathcal{F}}Q$ that satisfies
\[   T(Q/S)= \left\{\begin{array}{cc}
M(G/K)_{S/K} &  \mbox{if} \  K \subseteq S \\
0 & \mbox{otherwise} ,
\end{array}\right.\]
and such that there exists a surjection of Mackey functors
\[   \mathrm{coind}_Q^G(T) \rightarrow M .\]
\end{lemma}
\begin{proof} Given a finite subgroup $K$ of $G$ and its normalizer $Q$, there exists an inflation functor (see \cite[below Prop. 2.2]{ThevenazWebb95})
\[    \mathrm{inf}_{Q/K}^{Q}: \mathrm{Mack}_{\mathcal{F}}Q/K \rightarrow \mathrm{Mack}_{\mathcal{F}}Q : N \mapsto  \mathrm{inf}_{Q/K}^{Q}(N) \]
given by 
\[  \mathrm{inf}_{Q/K}^{Q}(N)(S)= \left\{  \begin{array}{cc} N(S/K) & \mbox{if $K \subseteq S$} \\ 0 & \mbox{otherwise.} \end{array}\right.  \]
Moreover, if $M \in \mbox{Mack}_{\mathcal{F}}G$ is a Mackey functor that is generated by $K$ such that
$M(G/S)=0$ for all $S \subsetneq K$, then by Section 4 in \cite{MartinezNucinkis06}, there is an adjunction isomorphism 
\[  \mathrm{Hom}_{Q/K}(U,M(G/K))\cong \mathrm{Hom}_{\mathcal{M}_{\mathcal{F}}}(\mathrm{ind}_Q^G\mathrm{inf}^Q_{Q/K}(\overline{U}),M)      \]
for every $Q/K$-module $U$.  Here $\overline{U}$ is the coinvariance functor associated to the $Q/K$-module $U$. Hence, we can choose $U=M(G/K)$ to obtain a map 
\[   \mathrm{ind}_Q^G\mathrm{inf}^Q_{Q/K}(\overline{M(G/K)}) \rightarrow M \]
Denoting $T=\mathrm{inf}^Q_{Q/K}(\overline{M(G/K)})$, we get a map
\[ \psi: \mathrm{ind}_Q^G(T) \rightarrow M. \]
By using the fact that $T(Q/K)=M(G/K)$, one can verify that $\psi$ is surjective when evaluated at $G/K$. Since $M$ is generated by $K$, we conclude that $\psi$ is surjective. It now suffices to prove that $\mathrm{ind}_Q^G(T)=\mathrm{coind}_Q^G(T)$. By Proposition \ref{prop: key mackey1}, this is equivalent to showing that the natural inclusion
\begin{equation}\label{eq: prod coprod}\bigoplus_{x \in [Q \setminus G/S]}T(Q/{}^xS\cap Q) \rightarrow  \prod_{x \in [Q \setminus G/S]}T(Q/{}^xS\cap Q)  \end{equation}
is an isomorphism for every finite subgroup $S$ of $G$. Now, by construction $T(Q/{}^xS\cap Q)$ is zero unless  $K \subseteq {}^xS\cap Q$. Since $S$ is a finite group, it can contain only finitely many different subgroups of the form $K^x$, with $x \in G$. On the other hand, if $K^x=K^y$, then $xy^{-1} \in Q$ and hence $x$ and $y$ represent the same element in $[Q \setminus G /S]$. We conclude that for every finite subgroup $S$ of $G$, there are only finitely many $x \in [Q\setminus G / S]$ such that $T(Q/{}^xS\cap Q)$ is non-zero. This shows that the map (\ref{eq: prod coprod}) is an isomorphism for every finite subgroup $S$ of $G$. 
\end{proof}
\begin{lemma} \label{th: weyl} Let $G$ be a group such that $\underline{\mathrm{cd}}_{\mathcal{M}}(G)=n<\infty$. Let $M \in \mbox{Mack}_{\mathcal{F}}G$ be generated by a collection of finite subgroups $\Omega$, such that there is a uniform bound on the lengths of the groups in $\Omega$. If $\mathrm{H}^n_{\mathcal{F}}(G,M) \neq 0$, then there exists a finite subgroup $F$ of $G$ such that $M(G/F)\neq 0$ and a $W_G(F)$-module $V$ such that
\[     \mathrm{H}^n_{\mathcal{F}}(W_G(F),\overline{V}) \neq 0.       \]
\end{lemma}

\begin{proof}  By assumption, we have a surjection
\[  \bigoplus_{F \in \Omega} M_{F} \rightarrow M \rightarrow 0,\]
where $M_F$ is the subfunctor of $M$ generated by $F \in \Omega$. A long exact cohomology sequence implies that $\mathrm{H}^{n}_{\mathcal{F}}(G,\bigoplus_{F \in \Omega} M_{F}^{\ast})\neq 0$. 
We will now make two reductions. \\

\noindent \textbf{A. We may assume that for each $F \in \Omega$ one has $M_F(G/S)=0$ for all $S \subsetneq F$.}
\begin{itemize}
\item[]  We claim that by possibly modifying the collection $\Omega$ and functors $M_F$, we may assume that for each $F \in \Omega$ one has $M_F(G/S)=0$ for all $S \subsetneq F$, while $\mathrm{H}^{n}_{\mathcal{F}}(G,\bigoplus_{F \in \Omega} M_{F}^{\ast}) \neq 0$ still holds. We will prove this claim by induction on 
\[ \tau(\{M_F\}_{F\in\Omega})= \max\{   l(F) \ | \ F \in \Omega \ \mbox{s.t.} \ \exists S \subsetneq F : M_F(G/S)\neq 0       \} . \]
If $M_F(G/S)=0$ for all $S \subsetneq F$ and for all $F \in \Omega$, then we set $ \tau(\{M_F\}_{F\in\Omega})=0$.
Note that $\tau(\{M_F\}_{F\in\Omega})$ is finite, since there is a uniform bound on the lengths of the groups in $\Omega$. If $\tau(\{M_F\}_{F\in\Omega})$ is zero, then there is nothing to prove. Let $k>0$ and assume our claim is true when $\tau(\{M_F\}_{F\in\Omega})<k$. Now proceed by induction and assume that $\tau(\{M_F\}_{F\in\Omega})=k$ and $\mathrm{H}^{n}_{\mathcal{F}}(G,\bigoplus_{F \in \Omega} M_{F}^{\ast}) \neq 0$.

Let $F \in \Omega$ and define $\tilde{\Omega}_F$ to be the set containing all $S \subsetneq F$ for which $M_F(G/S)\neq 0$. If no such $S$ exist then set $\tilde{\Omega}_F=\{F\}$. For each $S \in \tilde{\Omega}_F$, define $N_{S}$ to be the Mackey subfunctor of $M_F$ generated by $S$. Let $\tilde{\Omega}=\coprod_{F \in \Omega} \tilde{\Omega}_F$ and define $N_F=\bigoplus_{S \in \tilde{\Omega}_{F}}N_{S}$.
Finally, consider the natural map from $N_F$ to $M_F$ and denote the kernel of this map by $R_F$. Now consider the short exact sequences
\begin{equation} \label{eq: first one}   0 \rightarrow \bigoplus_{F \in \Omega} R_F \rightarrow \bigoplus_{F \in \Omega} N_F  \rightarrow  \bigoplus_{F \in \Omega} B_F \rightarrow 0        \end{equation}
and 
\begin{equation} \label{eq: second one}   0 \rightarrow \bigoplus_{F \in \Omega} B_F \rightarrow \bigoplus_{F \in \Omega} M_F  \rightarrow  \bigoplus_{F \in \Omega} M_F/B_F \rightarrow 0.        \end{equation}
The long exact cohomology sequence associated to (\ref{eq: second one}) implies that either
\[\mathrm{H}^{n}_{\mathcal{F}}(G,\bigoplus_{F \in \Omega} (M_F/B_{F})^{\ast}) \neq 0\]
or
\[\mathrm{H}^{n}_{\mathcal{F}}(G,\bigoplus_{F \in \Omega} B_{F}^{\ast}) \neq 0.\]

In the first case, we have $M_F/B_{F}(G/S)=0$ for all $S \subsetneq F$ and $M_F/B_{F}$ is generated by $F$, as desired.
In the second case we obtain that $\mathrm{H}^{n}_{\mathcal{F}}(G,\bigoplus_{F \in \Omega} N_{F}^{\ast}) \neq 0$, by considering the long exact cohomology sequence associated to (\ref{eq: first one}).

Since 
\[\bigoplus_{F \in \Omega} N_{F}= \bigoplus_{S \in \tilde{\Omega}} N_S ,  \]
and by construction $ \tau(\{N_S\}_{S\in\tilde{\Omega}})<k$, the claim follows by induction. 
\end{itemize}
\noindent \textbf{B. We may assume that $\Omega$ contains at most one representative of each conjugacy class of subgroups of $G$, while assumption A still holds.}
\begin{itemize}
\item[] To prove this, let $\tilde{\Omega}$ be a subset of $\Omega$ containing exactly one representative from each conjugacy class of subgroups in $\Omega$ and define \[\tilde{M}_{F}=\bigoplus_{\substack{F_0 \in \Omega \ \mbox{\tiny s.t.}\\ \exists g \in G: F_0=F^g}}=M_{F_0}\] for each $F \in \tilde{\Omega}$. Then clearly 
\[  \bigoplus_{F \in \tilde{\Omega}} \tilde{M}_{F}=  \bigoplus_{F \in \Omega} M_{F}.       \]
Moreover $\tilde{M}_{F}$ is generated by $F$ and $\tilde{M}_{F}(G/S)=0$ for every subgroup $S \subsetneq F$, since $D(G/S)\cong D(G/S^{g})$ for every Mackey functor $D$ and every $g \in G$. This finishes the proof of reduction $B$.
\end{itemize}
We now proceed assuming $A$ and $B$.

Denote, for each $F \in \Omega$, the normalizer of $F$ in $G$ by $Q_F$. By Lemma \ref{lemma: inflation} there exists for each $F \in \Omega$ a Mackey functor $T_F  \in  \mbox{Mack}_{\mathcal{F}}Q_F$ that satisfies
\[   T_F(Q_F/S)= \left\{\begin{array}{cc}
M_F(G/F)_{S/F} &  \mbox{if} \  F \subseteq S \\
0 & \mbox{otherwise,} 
\end{array}\right.\]
and a surjection of Mackey functors
\[   R_F=\mathrm{coind}_{Q_F}^G(T_F) \rightarrow M_F .\]
One easily checks, using Proposition \ref{prop: key mackey1}, that for each finite subgroup $S$ of $G$, the group $R_F(G/S)$ is non-trivial only if $F$ is conjugate to a subgroup of $S$. Since $\Omega$ contains at most one representative of each conjugacy class of subgroups of $G$ it follows that for each finite subgroup $S$ of $G$, $R_F(G/S)$ is non-zero for at most finitely many $F \in \Omega$. This shows that the natural inclusion
\[\bigoplus_{F \in \Omega}R_F \rightarrow \prod_{F\in \Omega} R_F \] is surjective.  Therefore, we obtain a surjection 
\[\prod_{F \in \Omega} R_F \rightarrow \bigoplus_{F \in \Omega} M_F\rightarrow 0. \]
A long exact cohomology sequence and the fact that cohomology commutes with products, now imply that there exists an $F \in \Omega$ such that $\mathrm{H}^{n}_{\mathcal{F}}(G,R_F^{\ast}) \neq 0 $. Since $R_F=\mathrm{coind}_{Q_F}^G(T_F)$ and $Q_F$ is the normalizer $N_G(F)$ of $F$ in $G$ , it follows from Shapiro's lemma and (\ref{eq: Mackey-bredon coind}) that $\mathrm{H}^{n}_{\mathcal{F}}(N_G(F),T_F^{\ast}) \neq 0 $.
Now consider the short exact sequence
\[   1 \rightarrow F \rightarrow N_G(F) \xrightarrow{\pi} W_G(F) \rightarrow 1.      \]
From \cite[Th. 5.1]{Martinez}, we obtain a convergent spectral sequence
\[   \mathrm{H}^p_{\mathcal{F}}(W_G(F),\mathrm{H}^q_{\mathcal{F} }( \pi^{-1}(-),T_F^{\ast})) \Longrightarrow \mathrm{H}^{p+q}_{\mathcal{F}}(N_G(F),T_F^{\ast}).        \]
Since $\pi^{-1}(K)$ is finite, for every finite subgroup $K$ of $W_G(F)$, this spectral sequence collapses and yields
\[  \mathrm{H}^n_{\mathcal{F}}(W_G(F),T_F(N_G(F)/\pi^{-1}(-))^{\ast}) = \mathrm{H}^{n}_{\mathcal{F}}(N_G(F),T_F^{\ast})\neq 0.  \]
Since $M_F$ is obtained form $M$ by taking a certain combination of successive subfunctors and quotients, and $T_F(N_G(F)/\pi^{-1}(-))=\overline{M_F(G/F)}$ by construction, the theorem is proven.
\end{proof}
We can now prove the second part of Theorem A.
\begin{theorem}\label{th: mackey dim eq} Let $G$ be a group of finite length, then 
\[    \mathcal{F}\mathrm{cd}(G) = \underline{\mathrm{cd}}_{co\mathcal{M}}(G) = \underline{\mathrm{cd}}_{\mathcal{M}}(G).\]
\end{theorem}
\begin{proof} We may assume that $ \mathcal{F}\mathrm{cd}(G)<\infty$, otherwise there is nothing to prove. Since $l(G)< \infty$, Theorem \ref{th: length orbit} implies that $\underline{\mathrm{cd}}_{co\mathcal{M}}(G) \leq \underline{\mathrm{cd}}_{\mathcal{M}}(G)\leq  \underline{\mathrm{cd}}(G)<\infty$. It follows from Proposition \ref{prop: Fcd2} that  $\mathcal{F}\mathrm{cd}(G) = \underline{\mathrm{cd}}_{co\mathcal{M}}(G) $. Finally, let $n=\underline{\mathrm{cd}}_{\mathcal{M}}(G)$.
By Lemma \ref{lemma: finite kernel}, we have $ \underline{\mathrm{cd}}_{co\mathcal{M}}(W_G(F)) = \underline{\mathrm{cd}}_{co\mathcal{M}}(N_G(F))$ and therefore $ \underline{\mathrm{cd}}_{co\mathcal{M}}(W_G(F)) \leq \underline{\mathrm{cd}}_{co\mathcal{M}}(G)$ for every finite subgroup $F$ of $G$. Since coinvariance functors are cohomological Mackey functors, Lemma \ref{th: weyl} implies that $\underline{\mathrm{cd}}_{co\mathcal{M}}(G) =\underline{\mathrm{cd}}_{\mathcal{M}}(G) $.\\

\end{proof}
Theorem A naturally leads to the following questions, which partially appear also in the open questions section of \cite{MartinezNucinkis06}. These questions can be seen as refinements of Nucinkis' conjecture mentioned in the introduction.
\begin{question} \rm Let $G$ be a group with $l(G)=\infty$.
\begin{itemize}
\item[(1)] Can one have $\mathcal{F}\mathrm{cd}(G)<\infty$ but $\underline{\mathrm{cd}}_{co\mathcal{M}}(G)=\infty?$ 
\item[(2)] Can one have $\underline{\mathrm{cd}}_{co\mathcal{M}}(G) < \underline{\mathrm{cd}}_{\mathcal{M}}(G)$ ?  Can one have $\underline{\mathrm{cd}}_{co\mathcal{M}}(G) < \infty$ but $\underline{\mathrm{cd}}_{\mathcal{M}}(G)=\infty$?
\item[(3)] Can one have $\underline{\mathrm{cd}}_{\mathcal{M}}(G) < \infty$ but $\underline{\mathrm{cd}}(G)=\infty$?
\end{itemize}

\end{question}
\section{Closure properties of $\mathfrak{F}$ and $\mathfrak{M}$}
We begin with the following lemma that deals with extensions with finite kernels. The statement about $\underline{\mathrm{cd}}(\Gamma)$ is well known, while the result about $\mathcal{F}\mathrm{cd}(G)$ also partially follows from \cite[Prop. 3.8]{Gandini12}.
\begin{lemma} \label{lemma: finite kernel} Consider a short exact sequence of groups
\[  1 \rightarrow F \rightarrow \Gamma \xrightarrow{\pi} G \rightarrow 1  \]
with $F$ a finite group. One has $\underline{\mathrm{cd}}(\Gamma)= \underline{\mathrm{cd}}(G)$, $\underline{\mathrm{cd}}_{(co)\mathcal{M}}(\Gamma)= \underline{\mathrm{cd}}_{(co)\mathcal{M}}(G)$, and  $\mathcal{F}\mathrm{cd}(G)=\mathcal{F}\mathrm{cd}(\Gamma)$.
\end{lemma}
\begin{proof}
By \cite[Th. 5.1]{Martinez}, we have a convergent spectral sequence
\[    \mathrm{E}_2^{p,q}(M)=\mathrm{H}^p_{\mathcal{F}}(G,\mathrm{H}^q_{\mathcal{F}}(\pi^{-1}(-),M)) \Longrightarrow \mathrm{H}^{p+q}_{\mathcal{F}}(\Gamma,M)       \]
for every right $\mathcal{O}_{\mathcal{F}}\Gamma$-module $M$. Since $\pi^{-1}(H)$ is finite for every finite subgroup $H$ of $G$, we have $\mathrm{H}^q_{\mathcal{F}}(\pi^{-1}(H),M)=0$ for all $q>0$ and $\mathrm{H}^0_{\mathcal{F}}(\pi^{-1}(H),M)=M(G/\pi^{-1}(H))$. This implies that the spectral sequence collapses and yields
\[   \mathrm{H}^n_{\mathcal{F}}(G,M(G/\pi^{-1}(-))) \cong \mathrm{H}^{n}_{\mathcal{F}}(\Gamma,M)    \]
for every right $\mathcal{O}_{\mathcal{F}}\Gamma$-module $M$. We conclude that $\underline{\mathrm{cd}}(\Gamma)\leq \underline{\mathrm{cd}}(G)$. Now let $N$ be a right $\mathcal{O}_{\mathcal{F}}G$-module. The map $\pi: \Gamma \rightarrow G$ induces a functor \[\pi: \mathcal{O}_{\mathcal{F}}\Gamma \rightarrow \mathcal{O}_{\mathcal{F}}G: \Gamma/H \mapsto G/\pi(H).\] Denote by $\mathrm{res}_{\pi}$ the associated restriction functor for right Bredon modules and let $M=\mathrm{res}_{\pi}(N)$. Then $M(G/\pi^{-1}(-))=N$ and hence
\[   \mathrm{H}^n_{\mathcal{F}}(G,N) \cong \mathrm{H}^{n}_{\mathcal{F}}(\Gamma,M)    \]
for every right $\mathcal{O}_{\mathcal{F}}G$-module $N$. It follows that $\underline{\mathrm{cd}}(G)\leq \underline{\mathrm{cd}}(\Gamma)$ and hence $\underline{\mathrm{cd}}(G)= \underline{\mathrm{cd}}(\Gamma)$. \\
\indent The other algebraic equalities are proven in a similar fashion as the equality for Bredon cohomological dimension, by restricting to (cohomological) Mackey functor coefficients or fixed point funtor coefficients.
\end{proof}

\begin{lemma} \label{lemma: dealing with N} Let $N$ be a finite index subgroup of $G$ and assume that $\underline{\mathrm{cd}}_{\mathcal{M}}(G)=n$. The following hold.
\begin{itemize}
\item[(1)]Let $M \in \mbox{Mack}_{\mathcal{F}}G$ be a Mackey functor generated by finite subgroups $F \subset N$.  If $\mathrm{H}^{n}_{\mathcal{F}}(G,M^{\ast})\neq 0$ then  $\mathrm{H}^{n}_{\mathcal{F}}(N,M^{\ast})\neq 0$ and hence $\underline{\mathrm{cd}}_{\mathcal{M}}(N)=n$.
\item[(2)] If $\mathrm{H}^{n}_{\mathcal{F}}(G,M^{\ast})= 0$ for all Mackey funtors $M \in \mbox{Mack}_{\mathcal{F}}G$ for which $M(G/F)=0$ for all $F \subseteq N$, then $\underline{\mathrm{cd}}_{\mathcal{M}}(N)=n$.
\end{itemize}
Moreover, (1) and (2) remain valid when restricting to cohomological Mackey functors and assuming that $\underline{\mathrm{cd}}_{co\mathcal{M}}(G)=n$.
\end{lemma}
\begin{proof} Let us first prove (1).
Let $\mathrm{ind}_N^G$ ($\mathrm{res}_N^G$) be the induction  (restriction) functor associated to $\mathcal{M}_{\mathcal{F}}N \rightarrow \mathcal{M}_{\mathcal{F}}G$.  Since $\mathrm{ind}_N^G$ is left adjoint to  $\mathrm{res}_N^G$, we obtain a canonical map
\[  p: A=\mathrm{ind}_N^G (\mathrm{res}_N^G(M) )\rightarrow M. \]
One can easily verify that $p(G/F)$ is surjective for all $F \subseteq N$. Since $M$ is generated by finite subgroups of $N$, it follows that $p$ is surjective. A long exact cohomology sequence now yields the exact sequence
\[   \mathrm{H}^{n}_{\mathcal{F}}(G,A^{\ast}) \rightarrow \mathrm{H}^{n}_{\mathcal{F}}(G,M^{\ast}) \rightarrow 0. \]
Hence, we have $\mathrm{H}^{n}_{\mathcal{F}}(G,A^{\ast}) \neq 0$. Now Shapiro's lemma and Proposition \ref{prop: key mackey1} imply that
\[    \mathrm{H}^{n}_{\mathcal{F}}(N,M^{\ast})= \mathrm{H}^{n}_{\mathcal{F}}(G,A^{\ast}) \neq 0.     \]
This proves part (1) of the lemma. \\
\indent Now let $M \in \mbox{Mack}_{\mathcal{F}}G$ such that $\mathrm{H}^n_{\mathcal{F}}(G,M^{\ast})\neq 0$, denote \[R=\bigoplus_{\substack{F \subseteq N\\ |F|< \infty}} M_F \]and consider the following exact sequence of Mackey functors
\[  0 \rightarrow R_0 \rightarrow R \rightarrow M \rightarrow M_0 \rightarrow 0, \]
where $R_0$ and $M_0$ are the kernel a cokernel of the natural map from $R$ to $M$.
Since  $M_0(G/F)=0$ for all $F \subseteq N$, it follows from the assumptions of $(2)$ and the long exact cohomology sequence associated to $0 \rightarrow R/R_0 \rightarrow M \rightarrow M_0 \rightarrow 0$ that $\mathrm{H}^{n}_{\mathcal{F}}(G,(R/R_0)^{\ast})\neq 0$. But then the long exact cohomology sequence associated to $0 \rightarrow R_0 \rightarrow R \rightarrow R/R_0 \rightarrow 0$ implies that $\mathrm{H}^{n}_{\mathcal{F}}(G,R^{\ast})\neq 0$. Since $R$ is generated by the finite subgroups of $N$, it follows from part (1) of the lemma that $\underline{\mathrm{cd}}_{\mathcal{M}}(N)=n$. This proves part (2). \\

Note that induction and restriction of cohomological Mackey functors produces cohomological Mackey functors. Moreover, quotient and subfunctors and direct sums of cohomological Mackey functors are also cohomological Mackey functors. Therefore, the preceding arguments remain valid in the context of cohomological Mackey functors and $\underline{\mathrm{cd}}_{co\mathcal{M}}$.
\end{proof}
We are now ready to prove Theorem B from the introduction, which gives a partial answer to Question 6.3 in \cite{MartinezNucinkis06}. 

\begin{proof}[proof of Theorem B] We shall only prove the statement for $\underline{\mathrm{cd}}_{\mathcal{M}}(G)$. The statement for $\underline{\mathrm{cd}}_{co\mathcal{M}}(G)$ is proven entirely similar, by restricting to cohomological Mackey functors. From \cite[Th. 5.1]{Martinez} we obtain a convergent spectral sequence
\[   \mathrm{H}^p_{\mathcal{F}}(Q,\mathrm{H}^q_{\mathcal{F}}(\pi^{-1}(-),M^{\ast})) \Longrightarrow \mathrm{H}^{p+q}_{\mathcal{F}}(G,M^{\ast}) \]
for every Mackey functor $M \in \mbox{Mack}_{\mathcal{F}}G$. It follows from Propostion \ref{prop: key mackey2} that
\[    \mathcal{M}_{\mathcal{F}}Q \rightarrow \mathrm{Ab}: Q/F \mapsto  \mathrm{H}^q_{\mathcal{F}}(\pi^{-1}(F),M^{\ast})  \]
is a Mackey functor. A standard spectral sequence argument now reduces the problem to the case where $Q$ is finite.

Assume $Q$ is a finite group. We will use induction on the index $[G:N]$. If $[G:N]=1$, then the statement is trivial. Now assume the statement is true for all extensions for which the order of the quotient group is strictly smaller than  $|Q|$.

We may assume that $N$ is the maximal non-trivial normal subgroup  of $G$, otherwise we would be done by induction. By assumption, we have $\underline{\mathrm{cd}}_{\mathcal{M}}(G)=n< \infty$ (if $\underline{\mathrm{cd}}(N)<\infty$ then $\underline{\mathrm{cd}}(G)<\infty$ by \cite[Th. 2.4]{Luck1} ) and let $M \in \mbox{Mack}_{\mathcal{F}}G$ be a Mackey functor for which $\mathrm{H}^{n}_{\mathcal{F}}(G,M^{\ast})\neq 0$. Our goal is to show that $n=\underline{\mathrm{cd}}_{\mathcal{M}}(N)$. We may proceed while assuming that $M(G/F)=0$ for all $F \subseteq N$, since otherwise we would be done by Lemma \ref{lemma: dealing with N}(2).

Choose a collection $\Omega$ of finite subgroups of $G$ such that $\Omega$ generates $M$. We can assume that $\Omega$ contains at most one representative of each conjugacy class of finite subgroups of $G$ and that $\Omega$ does not contain subgroups of $N$.  Since there is a uniform bound on the lengths of subgroups in $\Omega$, it follows from Lemma \ref{th: weyl}, that there exists a finite subgroup $F$ of $G$ that is not contained in $N$ such that $\underline{\mathrm{cd}}_{\mathcal{M}}(N_G(F))=\underline{\mathrm{cd}}_{\mathcal{M}}(W_G(F))=n$. We now consider two cases. First assume that $\pi(N_G(F))=Q$. Then $NF$ is a normal subgroup of $G$, which implies that $\pi(F)=Q$ since $N$ is maximal normal. It follows that \[W_G(H)=N_G(F)/F \cong N_G(F) \cap N\Big/ N\cap F.\]  and therefore 
\[n =\underline{\mathrm{cd}}_{\mathcal{M}}(W_G(F))= \underline{\mathrm{cd}}_{\mathcal{M}}(N_G(F) \cap N)\leq \underline{\mathrm{cd}}_{\mathcal{M}}( N).\] We conclude $\underline{\mathrm{cd}}_{\mathcal{M}}(N)=n.$

Secondly, assume that $\pi(N_G(F))=Q_0 \subsetneq Q$. In this case, there is a short exact sequence
\[    1 \rightarrow N_G(F)\cap N \rightarrow N_G(H) \rightarrow Q_0 \rightarrow 1   \]
By induction we conclude that $\underline{\mathrm{cd}}_{\mathcal{M}}(N_G(F)\cap N)=n$ and hence $\underline{\mathrm{cd}}_{\mathcal{M}}(N)=n$ which completes the induction and finishes the proof. 

\end{proof}

\begin{remark} \rm Let $G$ be a group with a finite index subgroup $N$. If $\underline{\mathrm{cd}}(N)< \infty$, then a geometric construction of L\"{u}ck (see \cite[Th. 2.4]{Luck1}) implies that $\underline{\mathrm{cd}}(G)\leq [G:N]\underline{\mathrm{cd}}(N)$. In particular, it follows that $\underline{\mathrm{cd}}(G)<\infty$. Such a result is missing for Mackey functors, meaning that we do not know whether $ \underline{\mathrm{cd}}_{(co)\mathcal{M}}(N)< \infty$ implies  $\underline{\mathrm{cd}}_{(co)\mathcal{M}}(G)< \infty$, if we do not assume that  $\underline{\mathrm{cd}}(N)< \infty$. This is the reason why we must assume that either $\underline{\mathrm{cd}}(N)<\infty$ or $\underline{\mathrm{cd}}_{(co)\mathcal{M}}(G)< \infty$ in the statement of Theorem B. On other hand, Nucinkis' conjecture implies that $\underline{\mathrm{cd}}_{(co)\mathcal{M}}(N)< \infty$ if and only if $\underline{\mathrm{cd}}(N)< \infty$.

\end{remark}
Therefore, the following question is of some interest. 

\begin{question}\rm Let $G$ be a group with a finite index subgroup $N$. Can one prove that
\[\underline{\mathrm{cd}}_{(co)\mathcal{M}}(N) <\infty \Longrightarrow   \underline{\mathrm{cd}}_{(co)\mathcal{M}}(G)<\infty,    \]
without assuming that $\underline{\mathrm{cd}}(N)<\infty$?

\end{question}
Finally, we turn to the proof of Theorem C which requires the following lemma.
\begin{lemma} \label{lemma: spec}A contractible finite dimensional $G$-CW-complex $X$ and a $G$-module $M$ give rise to a convergent $E_1$-term spectral sequence
\[    \prod_{\sigma \in \Delta_p} \mathrm{H}^q_{\mathcal{F}}(G_{\sigma},\underline{M}^{\ast})    \Longrightarrow   \mathrm{H}^{p+q}_{\mathcal{F}}(G,\underline{M}^{\ast}) ,   \]
 where $\Delta_p$ contains representatives of $G$-orbits of $p$-dimensional cells of $X$ and $G_{\sigma}$ is the stabilizer of the cell $\sigma$. 
\end{lemma}
\begin{proof}
In the proof of Proposition $3.2$ of \cite{DemPetTal} , a convergent $E_1$-spectral sequence
\begin{equation*} \label{eq: spec}    \prod_{\sigma \in \Delta_p} \mathrm{H}^q_{\mathcal{F}}(G_{\sigma},V)    \Longrightarrow  \mathrm{H}^{p+q}(\mathrm{Hom}_{\mathcal{F}}(C_{\ast}(X^{-}) \otimes C_{\ast}(\underline{E}G^{-})),V))      \end{equation*}
is constructed for every right $\orb$-module $V$. Here $C_{\ast}((X \times \underline{E}G)^H )=C_{\ast}(X^H) \otimes C_{\ast}(\underline{E}G^H)$ is the cellular chain complex of the $H$-fixed point set $(X \times \underline{E}G)^H$, $\Delta_p$ contains representatives of $G$-orbits of $p$-dimensional cells of $X$ and $G_{\sigma}$ is the stabilizers of the cell $\sigma$. 
We will apply this for $V=\underline{M}^{\ast}$. Clearly one has 
\[     \mathrm{H}^{r}(\mathrm{Hom}_{\mathcal{F}}(C_{\ast}(X^{-}) \otimes C_{\ast}(\underline{E}G^{-}),\underline{M}^{\ast}))  =   \mathrm{H}^{r}(\mathrm{Hom}_{G}(C_{\ast}(X) \otimes C_{\ast}(\underline{E}G),M ))    \]
by adjointness of fixed point functors and evalution at the identity functors. We claim that
\[  \mathrm{H}^{r}_{\mathcal{F}}(G,\underline{M}^{\ast}) \cong \mathrm{H}^{r}(\mathrm{Hom}_{G}(C_{\ast}(X) \otimes C_{\ast}(\underline{E}G),M )           \]
for every $r \in \mathbb{N}$. To prove this, note that 
$\mathrm{H}^{r}_{\mathcal{F}}(G,\underline{M}^{\ast})=\mathcal{F}\mathrm{H}^{r}(G,M)$ by Proposition \ref{prop: rel coh}. Hence, since $X \times \underline{E}G$ is a contractible proper $G$-CW-complex, it suffices to show that $C_{\ast}(X) \otimes C_{\ast}(\underline{E}G)\rightarrow \mathbb{Z}$ is $H$-split when restricted to a finite subgroup $H$ of $G$. By Corollary 3.4 in \cite{nucinkis00}, we know that  $C_{\ast}(\underline{E}G) \rightarrow \mathbb{Z}$ is $H$-split, and by Theorem 1.4 in \cite{KrophollerWall11} the complex $C_{\ast}(X) \rightarrow \mathbb{Z}$ is also $H$-split. We therefore conclude that $C_{\ast}(X) \otimes C_{\ast}(\underline{E}G) \rightarrow \mathbb{Z}$ is $H$-split by Lemma 2.3 in \cite{nucinkis00}. This proves the lemma.

\end{proof}

\begin{proof}[proof of Theorem C]
Property (1) for $\mathcal{F}_d$ is proven in \cite[Lemma 2.1]{nucinkis00}, and also follows from Shapiro's lemma, (\ref{eq: Mackey-bredon coind}) and Lemma \ref{lemma: ind fixed point}. Property (1) for $\mathfrak{M}_d$ is proven in (\ref{eq: mackey subgroups}), while Property (2) follows from Lemma \ref{lemma: finite kernel}. Property (3) is immediate from Theorem B, since all groups in $\mathfrak{L}\cap \mathfrak{M}_c$ have finite Bredon cohomological dimension for proper actions, by Theorem A.

Now assume that $X$ is a contractible $G$-CW-complex of dimension at most $n$ and with isotropy groups in $\mathfrak{F}_d$. For a $G$-module $M$, it follows form the spectral sequence of the Lemma \ref{lemma: spec} that $\mathrm{H}^{n+d+1}_{\mathcal{F}}(G,\underline{M}^{\ast})=0$. This proves Property (4).

Finally, assume that $G$ that act by isometries, discrete orbits and with isotropy in $\mathfrak{F}_d$ or $\mathfrak{M}_d$ on a separable CAT($0$)-space $Y$ of topological dimension at most $n$.
Following the proof of Theorem A in \cite{DP3} one can use Proposition 2.6 in \cite{DP3} to find a $G$-CW-complex $X$ of dimension at most $n$ and with isotropy in $\mathfrak{F}_d$ or $\mathfrak{M}_d$ together with cellular $G$-maps
\[      \underline{E}G \xrightarrow{f} X \times \underline{E}G \xrightarrow{g} \underline{E}G           \]
such that $g \circ f$ is $G$-homotopy equivalent to the identity map via cellular maps. By passing to cellular $\orb$-chain complexes, we obtain chain maps
\[      C_{\ast}(\underline{E}G^{-}) \xrightarrow{f} C_{\ast}(X^{-}) \otimes C_{\ast}(\underline{E}G^{-}) \xrightarrow{g} C_{\ast}(\underline{E}G^{-})           \]
such that $g \circ f$ is chain homotopy equivalent to the identity map. This implies that for every right $\orb$-module $M$ and every $r \in \mathbb{N}$, we have maps
\[     \mathrm{H}^{r}_{\mathcal{F}}(G,M)=\mathrm{H}^r(\mathrm{Hom}_{\mathcal{F}}(C_{\ast}(\underline{E}G^{-}),M)) \rightarrow  \mathrm{H}^r(\mathrm{Hom}_{\mathcal{F}}(C_{\ast}(X^{-})\otimes C_{\ast}(\underline{E}G^{-}),M))   \rightarrow    \mathrm{H}^{r}_{\mathcal{F}}(G,M)        \]
whose composition is the identity map. We deduce that $  \mathrm{H}^{r}_{\mathcal{F}}(G,M)$ injects into \[ \mathrm{H}^r(\mathrm{Hom}_{\mathcal{F}}(C_{\ast}(X^{-})\otimes C_{\ast}(\underline{E}G^{-}),M)) \] for every $r \in \mathbb{N}$. As noted in the proof of Lemma \ref{lemma: spec}, there is a convergent $E_1$-spectral sequence
\[    \prod_{\sigma \in \Delta_p} \mathrm{H}^q_{\mathcal{F}}(G_{\sigma},M)    \Longrightarrow  \mathrm{H}^{p+q}(\mathrm{Hom}_{\mathcal{F}}(C_{\ast}(Y^{-}) \otimes C_{\ast}(\underline{E}G^{-})),M)      \]
for every $M \in \orb\mbox{-Mod}$. Property $(5)$, respectively $(6)$, now follows from this spectral sequence by restricting to fixed point functor coefficients, respectively Mackey functor coefficients, and assuming that $Y$ (and hence $X$) has isotropy in $\mathfrak{F}_d$, respectively $\mathfrak{M}_d$.

\end{proof}

The following question now also presents itself.

\begin{question} Does property (3) of Theorem $C$ remain valid when we remove $\mathfrak{L}$? Does property (4) of Theorem C remain valid if we replace $\mathfrak{F}$ with $\mathfrak{M}$?

\end{question}

We finish this paper by pointing some consequences of the conjecture by Nucinkis stating that
\begin{equation} \label{eq: nucon}  \mathcal{F}\mathrm{cd}(G)<\infty \Longrightarrow \underline{\mathrm{cd}}(G)<\infty          \end{equation}
for all groups $G$. For a countable group $N$ with $\underline{\mathrm{cd}}(N)<\infty$, let $\mathfrak{E}_N$ be the class of all groups that arise as extensions of subgroups of $N$ by finite groups, and define the invariant 
\[   \delta(N)=\sup\{  \underline{\mathrm{cd}}(\Gamma) \ | \ \Gamma \in \mathfrak{E}_N  \} .     \]
A more restrictive version  of this invariant (that only allows extensions by finite cyclic groups) was defined in \cite{DP2} as a tool for studying the behavior under group extensions of Bredon cohomological dimension for the family of virtually cyclic subgroups (denoted by $\underline{\underline{\mathrm{cd}}}$). There are no groups $N$ known for which $\underline{\mathrm{cd}}(N)<\infty$ but $\delta(N)=\infty$, but it is in general a very difficult problem to determine whether or not $\delta(N)$ is finite for a given group $N$.  On the other hand, let $N \in \mathfrak{L}\cap \mathfrak{F}$ and assume that $\delta(N)=\infty$. Then there must exist a collection of groups $\{\Gamma_i \}_{i \in \mathbb{N}}$ such that $\Gamma_i$ has a finite index subgroup contained in $N$ and $\underline{\mathrm{cd}}(\Gamma_i)\geq i$. It follows property (3) of Theorem C that $\mathcal{F}\mathrm{cd}(\Gamma_i)\leq \mathcal{F}\mathrm{cd}(N)<\infty$ for every $i \in \mathbb{N}$. Now define $G$ as the infinite free product of all the groups $\Gamma_i$. Then clearly one has $\underline{\mathrm{cd}}(G)=\infty$. On the other hand, by considering the action of $G$ on the Bass-Serre tree associated to this free product, one obtains form property (4) of Theorem $C$ that $\mathcal{F}\mathrm{cd}(G)\leq 1+\mathcal{F}\mathrm{cd}(N)<\infty$. We are led to the following observation.
\begin{observation} \rm Assuming that (\ref{eq: nucon}) is true, the following three statements hold.

\begin{itemize}
\item[(1)] One has $\delta(N)<\infty$ for all $N \in \mathfrak{L}\cap\mathfrak{F}$.
\item[]
\item[(2)] The class of groups containing all groups with finite length and finite Bredon cohomological dimension for the family of finite subgroups, is closed under group extensions.
\item[]
\item[(3)] The class of group containing all groups with a uniform bound on the orders of their finite subgroups and with finite Bredon cohomological dimension for the family of virtually cyclic subgroups, is closed under group extensions.
\end{itemize}
Observation (1) follows from the discussion above. Observation (2) follows from observation (1) and a simple spectral sequence argument (see \cite[Cor 5.2]{Martinez}). Finally, observation (3) follows from combining observation (1), \cite[Th. B]{DP2} and \cite[Th. 2.4]{Luck1} together with a similar spectral sequence argument as before.
\end{observation}

\end{document}